\documentclass[11pt]{article}%
\usepackage{makeidx}
\usepackage{amssymb}
\usepackage[all]{xy}
\usepackage[latin1]{inputenc}
\usepackage{amsfonts}
\usepackage{amsmath}
\usepackage{amssymb}
\usepackage{url}
\usepackage{hyperref}
\usepackage{graphicx}%
\usepackage{color}
\usepackage{xcolor}
\usepackage{amsbsy}
\usepackage[all]{xy}
\usepackage{makeidx}
\usepackage{color}
\usepackage{mathrsfs}

\usepackage{color}

\providecommand{\U}[1]{\protect\rule{.1in}{.1in}}

\numberwithin{equation}{section}
\providecommand{\U}[1]{\protect\rule{.1in}{.1in}}
\providecommand{\U}[1]{\protect\rule{.1in}{.1in}}
\newtheorem{theo}{Theorem}[section]
\newtheorem{prop}[theo]{Proposition}
\newtheorem{lem}[theo]{Lemma}

\newtheorem{rem}[theo]{Remark}

\newcommand{\CC}{\mathbb{C}}

\newcommand{\EE}{\mathbb{E}}
\newcommand{\FF}{\mathbb{F}}
\newcommand{\HH}{\mathbb{H}}
\newcommand{\KK}{\mathbb{K}}
\newcommand{\LL}{\mathbb{L}}

\newcommand{\NN}{\mathbb{N}}

\newcommand{\Ba}{ {\cal B }}
\newcommand{\Ca}{ {\cal C }}
\newcommand{\Da}{ {\cal D }}

\newcommand{\Ea}{ {\cal E }}

\newcommand{\Fa}{ {\cal F }}

\newcommand{\Ga}{ {\cal G }}
\newcommand{\Qa}{ {\cal Q }}
\newcommand{\Ia}{ {\cal I }}
\newcommand{\Xa}{ {\cal X }}
\newcommand{\Ma}{ {\cal M }}

\newcommand{\Sa}{ {\cal S}}
\newcommand{\Ha}{ {\cal H }}

\newcommand{\Pa}{ {\cal P }}
\newcommand{\Za}{ {\cal Z }}

\newcommand{\point}{\mbox{\LARGE .}}
\newcommand{\proof}{\noindent\mbox{\bf Proof:}\\}
\newcommand{\cqfd}{\hfill\blbx \\}
\def\blbx{\hbox{\vrule height 5pt width 5pt depth 0pt}\medskip}

\def \EE{\mathbb{E}}

\def \CC{\mathbb{C}}
\def \LL{\mathbb{L}}

\def \XX{\mathbb{X}}

\usepackage{bbm}

\def \KK{\mathbbm{K}}

\let\oldchi\chi
\renewcommand{\chi}[1][1.3pt]{%
  \mathrel{\raisebox{#1}{\scalebox{1.2}{$\oldchi$}}}%
}

\def \XX{\mathbbm{X}}

\newcounter{hypA}
\newenvironment{hypA}{\refstepcounter{hypA}\begin{itemize}
  \item[({\bf A\arabic{hypA}})]}{\end{itemize}}

\begin{document}

%\title{A Sharp First Order Analysis of Feynman-Kac Particle Models}

\begin{center}

{\Large \textbf{A Sharp First Order Analysis of Feynman-Kac Particle Models}}

\bigskip

BY HOCK PENG CHAN$^{1}$, PIERRE DEL MORAL$^{2}$, \& AJAY JASRA$^{1}$

{\footnotesize $^{1}$Department of Statistics \& Applied Probability,
National University of Singapore, Singapore, 117546, SG.}\\
{\footnotesize E-Mail:\,}\texttt{\emph{\footnotesize stachp@nus.edu.sg, staja@nus.edu.sg}}\\
{\footnotesize $^{2}$School of Mathematics and Statistics,
University of New South Wales, Sydney NSW, 2052, AUS.}\\
{\footnotesize E-Mail:\,}\texttt{\emph{\footnotesize p.del-moral@unsw.edu.au}}
\end{center}

%\author{PIERRE DEL MORAL,\\ 
%\footnotesize{School of Mathematics and Statistics, University of New South Wales, \texttt{p.del-moral@unsw.edu.au}},\\ 
%AJAY JASRA,\\ 
%\footnotesize{Department of Statistics and Applied Probability, National University of Singapore, \texttt{staja@nus.edu.sg}}
%}
%\maketitle

\begin{abstract}
This article provides a new theory for the analysis of forward and backward particle approximations of Feynman-Kac models.
Such formulae are found in a wide variety of applications and their numerical (particle) approximation are required due
to their intractability. Under mild assumptions, we provide sharp and non-asymptotic first order expansions of these particle methods,
 potentially on path space and for possibly unbounded functions. 
 These expansions allows one to consider upper and lower bound bias type estimates for a given time horizon $n$ and particle number $N$; these non-asymptotic estimates are of order $\mathcal{O}(n/N)$.
Our approach is extended to tensor products of particle density profiles, leading to new sharp and non-asymptotic propagation of chaos 
estimates. The resulting upper and lower bound propagation of chaos estimates seems to be the first result of this kind for mean field particle models.
As a by-product of our results, we also provide some analysis of the particle Gibbs sampler, providing first order expansions of the kernel and minorization estimates.
\end{abstract}
%\tableofcontents

\section{Introduction}

Feynman-Kac formulae provide a very general description of several models, used in statistics, physics and many more; see \cite{d-2004,d-2013}.
For instance, in the context of non-linear filtering, they provide a precise characterization of the sequence of filtering and smoothing distributions. There are numerous other
interpretations of Feynman-Kac formulae, which can be found in many application areas as listed above, but essentially leads to sequence of probability measures. In most practical problems of interest, one cannot analytically compute the expressions associated to Feynman-Kac formula, nor the associated expectations. As a result, particle approximations
of such measures have been developed. 

Particle approximation methods use $N\geq 1$ samples (or particles) that are generated in parallel, sequentially in time, and are propagated via 
Markov kernels with subsequent selection (or resampling) steps; resulting in $n$ time steps of $N$ samples which are serially and with each-other dependent.
Extensions to the backward interpretation (e.g.~\cite{dds1}) have also been constructed; this latter interpretation is, for example, useful for the smoothing problem in non-linear filtering.
Several convergence results, as $N$ grows, have been proved; see for instance \cite{cl-2013,d-2004,d-2013,douc1}. In particular, as summarized in \cite{d-2004,d-2013}, a variety of results have been established, including the bias of the particle approximation of the marginal Feynman-Kac formula and associated propagation of chaos estimates (independence estimates
of $1\leq q<N$ particles). Some results for the backward interpretation can also be found in \cite{dds1,douc}.

In this article, under assumptions, we provide a first order expansion of the forward and backward particle approximation,
potentially on path space (i.e.~on the entire collection of time steps $0,\dots,n$). This expansion allows one to consider bias estimates (both upper and lower bounds) for possibly unbounded functions; these estimates are $\mathcal{O}(n/N)$.
The expansion is extended for tensor products of the particle approximation,  potentially on path space and for possibly unbounded functions, 
and new propagation of chaos properties are proved, again of $\mathcal{O}(n/N)$ and both upper and lower bounds.
The expansions, to our knowledge, are entirely new.
In addition, previous bias estimates did not hold for unbounded functions and were not established for the backward interpretation (all possibly on path spaces). 
Propagation of chaos properties were also not established for unbounded functions. In both cases, we are not aware of any lower-bounds in the literature.

As a by-product of our analysis, we consider the popular particle Gibbs (PG) sampler in \cite{adh-2010}. This is a Markov chain Monte Carlo algorithm designed to sample
from a single pre-specified probability distribution, associated to a single (particle approximation or backward particle approximation of a) Feynman-Kac formula (see \cite{adh-2010,dkp-14}) of $n$ time steps and given $N$. We provide a first order expansion of the PG Markov kernel, around its marginal target measure. This allows one to consider minorization estimates and rates
of convergence of the PG kernel; our estimates have been established elsewhere in \cite{alv-2013,dkp-14,ldm-2014} under different technical approaches and assumptions. Further discussion and comparison is given in section \ref{sec:dual}. In addition, first order propagation of chaos estimates are derived 
for empirical measures of the dual particle model with a frozen path; explanations are given in section \ref{sec:dual}.

We emphasize that the work in this article provides a new theory for the analysis of  forward and backward particle approximations of Feynman-Kac models. This is based upon our 
first order approach, which in turn, allows one to derive non-asymptotic and sharp propagation of chaos estimates; as noted above, these results did not previously exist in the literature.
In addition, the results have important practical implications, such as the bias estimate for the backward interpretation; the fact that the result holds on path space is very useful
for the smoothing problem in the non-linear filtering literature.
The proofs of our results use perturbation semigroup techniques  as well as empirical and combinatorial tensor product analysis.
This is new in comparison to much of the analysis that is summarized in \cite{d-2004,d-2013}.

This article is structured as follows. In section \ref{sec:fk_int} the notations of the paper and Feynman-Kac models are described.
In section \ref{sec:mean_field} the mean-field particle approximation is discussed and our main result (Theorem \ref{T2}) is stated. Some further developments
on propagation of chaos are also given. In section \ref{sec:dual}, the implication of our results for the particle Gibbs sampler are stated (Theorem \ref{T1}). In
section \ref{sec:semi} some original semigroup analysis is given, along with the proofs of Theorem \ref{T2} and \ref{T1}. Some technical results are given in the appendix
along with the proofs for our propagation of chaos results.

\section{Feynman-Kac Path Integration Models}\label{sec:fk_int}

\subsection{Notation}

Given some measurable space $S$ we denote respectively by $\Ma(S)$,
$\Pa(S)$ and $\Ba(S)$, the set of finite signed measures on $S$, the convex subset of probability measures, and the Banach space of bounded measurable functions equipped with the uniform norm $\Vert f\Vert=\sup_{\textsl{x}\in S}\vert f(\textsl{x})\vert$.  

The total variation norm on measures $\mu\in \Ma(S)$ 
 is defined by $$
 \Vert\mu\Vert_{\tiny tv}:=\sup_{f\in\Ba(S)~:~\|f\|\leq 1}\vert \mu(f)\vert
\quad\mbox{\rm with the Lebesgue integral}\quad
\mu(f):=\int~\mu(d\textsl{x})~f(\textsl{x}).
$$
 We also denote by  $\delta_a$ the Dirac measure at some state $a$, so that $\delta_a(f)=f(a)$. We say that $\nu\leq \mu$ as soon
 as $\nu(f)\leq \mu(f)$ for any non-negative function $f$.

A bounded integral operator $Q(x,dy)$ between the measurable spaces $S$ and $S^{\prime}$ is defined for any $f\in\Ba(S^{\prime})$
by the measurable function $Q(f)\in\Ba(S)$ defined by
$$
Q(f)(\textsl{x}):=\int~Q(\textsl{x},d\textsl{y})~f(\textsl{y}).
$$
The operator $Q$ generates a dual operator $\mu\in \Ma(S)\mapsto \mu Q\in\Ma(S^{\prime})$ by the dual formula
$
(\mu Q)(f)=\mu(Q(f))
$. 

When a bounded integral operator  $M$ from a state space $S$ into a possibly different state space $S^{\prime}$ has a constant mass, that is, when $M(1)\left(  x\right)  =M(1)\left(
y\right)  $ for any $(x,y)\in S^{2}$, the operator $\mu\mapsto\mu M$ maps the set
$\mathcal{M}_{0}(S)$ of measures $\mu$ on $S$ with null mass $\mu(1)=0$ into $\mathcal{M}_{0}(S^{\prime})$. In this situation, we let
$\beta(M)$ be the Dobrushin coefficient of a bounded integral operator $M$
defined by the formula
\begin{equation}\label{f27}
\beta(M):=\sup{\ \{\mbox{\rm osc}(M(f))\;;\;\;f~\mbox{\rm s.t.}~\mbox{\rm osc}(f)\leq 1\}} 
\end{equation} where $\mbox{\rm osc}(f):=\sup_{x,y}|f(x)-f(y)|$ stands for the oscillation of some
function.

The $q$-tensor product of $Q$ is the integral operator defined for any $f\in \Ba(S^q)$ by
$$
Q^{\otimes q}(f)(\textsl{x}^1,\ldots,\textsl{x}^q):=\int~\left\{\prod_{1\leq i\leq q}Q(\textsl{x}^i,d\textsl{y}^i)\right\}~f(\textsl{y}^1,\ldots,\textsl{y}^q).
$$
We also denote by $Q_1Q_2$ the composition of two operators defined by
$$
(Q_1Q_2)(\textsl{x},d\textsl{z}):=\int Q_1(\textsl{x},d\textsl{y})Q_2(\textsl{y},d\textsl{z}).
$$

The Boltzmann-Gibbs transformation $\Psi_G~:~\eta\in\Pa(S)\mapsto \Psi_{G}(\eta)\in \Pa(S)$ 
associated with some positive function $G$ on some state space $S$ is defined by
$$
\Psi_{G}(\eta)(d\textsl{x}):=\frac{1}{\eta(G)}~G(\textsl{x})~\eta(d\textsl{x}).
$$
We also denote by $\#(E)$ the cardinality of a finite set
and we use the standard conventions $\left(
\sup_{\emptyset},\inf_{\emptyset}\right)=\left(-\infty,+\infty\right)$, and $\left(\sum_{\emptyset},\prod_{\emptyset}\right)=(0,1)$.

\subsection{Feynman-Kac Models}

We consider a collection of bounded and positive potential functions
$G_n$ on some  measurable state spaces $S_n$, with $n\in\NN$.
We also let $X_n$ be a Markov chain on $S_n$ with initial distribution $\eta_0\in\Pa(S_0)$ and
some Markov transitions $M_n$ from $S_{n-1}$ into $S_n$. The Feynman-Kac measures $(\eta_n,\gamma_n)$  associated with the parameters
$(G_n,M_n)$ are defined for any $f_n\in\Ba(S_n)$ by
\begin{equation}\label{FK-def-intro-ref}
\eta_n(f_n):={\gamma_n(f_n)}/{\gamma_n(1)},\quad%\mbox{\rm with}\quad
\gamma_n(f_n):=\EE\left(f_n(X_n)~Z_n(X)\right), \quad%\mbox{\rm and}\quad 
Z_n(X):=\prod_{0\leq p<n}G_p(X_p).
\end{equation} 

The evolution equations associated with these measures are given by
\begin{equation}\label{evolve-eq-intro-ref}
\gamma_{n+1}=\gamma_{n}Q_{n+1}\quad\mbox{and}\quad
\eta_{n+1}=\Phi_{n+1}(\eta_n):=\Psi_{G_{n}}(\eta_{n})M_{n+1}
\end{equation} 
with the integral operators $$
Q_{n+1}(\textsl{x}_{n},d\textsl{x}_{n+1})=G_{n}(\textsl{x}_{n})~
M_{n+1}(\textsl{x}_{n},d\textsl{x}_{n+1}).
$$

The unnormalized measures $\gamma_n$ can be expressed in terms of the normalized ones
using the well known product formula
$$
\gamma_n(f_n)=\eta_n(f_n)~\prod_{0\leq p<n}\eta_p(G_p).
$$
We also recall the semigroup decompositions
$$
\forall 0\leq p\leq n\qquad \gamma_n=\gamma_pQ_{p,n}\qquad\mbox{\rm and}\quad \eta_n=\eta_p\overline{Q}_{p,n}
$$
with the integral operators $Q_{p,n}=Q_{p+1}\ldots Q_n$, and the normalized semigroups
$$
\overline{Q}_{p,n}(f_n)(\textsl{x}_p)={{Q}_{p,n}(f_n)(\textsl{x}_p)}/{\eta_p{Q}_{p,n}(1)}=(\overline{Q}_{p+1}\ldots \overline{Q}_n)(f_n)(\textsl{x}_p).
$$
In the above display, $\overline{Q}_{p+1}$ stands for the collection of integral operators defined as $Q_{p+1}$ by replacing $G_{p}$
with the normalized potential functions $\overline{G}_p=G_p/\eta_p(G_p)$.

Suppose that $X_n=(X_0^{\prime},\ldots,X_n^{\prime})$ is the historical process of an auxiliary 
Markov chain $X_n^{\prime}$ evolving in some measurable state space $S^{\prime}_n$ with some Markov transitions kernels $M^{\prime}_n$.
Also suppose that the potential functions $G^{\prime}_n(X^{\prime}_n)=G_n(X_n)$ only depends on the terminal state $X_n^{\prime}$ of the trajectory $X_n$.
In this context, the Feynman-Kac model (\ref{FK-def-intro-ref}) is called the historical version of the Feynman-Kac model associated with the parameters $(G_n^{\prime},M_n^{\prime})$.

\subsection{Illustrations}

Feynman--Kac models appear in numerous fields including signal processing, statistics, mathematical finance, biology, rare event analysis, chemistry and statistical physics; see
 \cite{cappe-moulines}, \cite{cdho}, \cite{d-2004}, \cite{d-2013}, and \cite{doucet2001}. Their interpretation 
depends on the application domain. For instance, in molecular and quantum physics, Feynman--Kac path integrals provide a probabilistic interpretation of
imaginary time Schr\"odinger equations. To be more precise, let $M^{\prime}_n\simeq_{\Delta t\downarrow 0}
Id+L~\Delta t$ be the Markov kernel  associated with some discretization of continuous-time stochastic
process $X^{\prime}_{t}$ with infinitesimal generator $L^{\prime}$ on some time mesh $t_{n+1}-t_n:=\Delta t\ll1$, with $t_n=n\lfloor t/n\rfloor $.
We also assume the existence of a potential $V$ and introduce the functions
$G^{\prime}_{n}=e^{-V~\Delta t}$. 

Replacing the chain $X^{\prime}_n$ above by the discrete time approximation model $X^{\prime}_{t_n}$, we have
$$
\eta_{n}(f)\propto\mathbb{E}\left(  f((X^{\prime}_{t_k})_{0\leq k\leq n})~e^{ -\sum_{0\leq t_k<t_n}V(X^{\prime}_{t_k}) (t_{k+1}-t_k) }\right)\simeq_{\Delta t\downarrow 0}
$$
\begin{equation}\label{FK-ex-qmc}
\mathbb{E}\left(  f((X^{\prime}_{s})_{ s\leq t})~e^{-\int_{0}^{t}%
V(X^{\prime}_{s})ds}\right).
\end{equation}
The marginal $\gamma_t$ w.r.t.~the terminal time $t$ of the r.h.s.~measures
is often defined, in a weak sense, by
the imaginary time Schr\"odinger equation $
\frac{d}{dt}\gamma_{t}(f)=\gamma_{t}(L^{V}(f))$, with $L^V(f)=L(f)-Vf$.
In computational physics and chemistry, the genetic particle models discussed above belong to the class of Quantum Monte Carlo methods. They are often termed Resampled Monte Carlo 
methods, or Diffusion Monte Carlo methodologies.
 For a more thorough
discussion of these continuous-time models and their applications in chemistry
and physics, see~\cite{caffarel,cances-tony,djlmp-2013,tony3,tony,tony-rs,rousset}, the recent monograph~\cite{d-2013}, and the
references therein.

\section{Mean Field Particle Models}\label{sec:mean_field}
\subsection{Description of the Models}
The mean field particle interpretation of the measures $(\eta_n,\gamma_n)$
starts with $N$ independent random variables $\xi_0:=(\xi^{i}_0)_{1\leq i\leq N}\in S_0^N$ with common law $\eta_0$.
The simplest way to evolve the population of $N$ individual (a.k.a.~samples, particles, or walkers)  $\xi_n:=(\xi^{i}_n)_{1\leq i\leq N}\in S_n^N$ is to consider $N$ conditionally independent individuals $\xi_{n+1}:=(\xi^{i}_{n+1})_{1\leq i\leq N}\in S_{n+1}^N$ with common distribution 
\begin{equation}\label{mean-field-intro}
\Phi_{n+1}(m(\xi_n)
)\quad\mbox{\rm with}\quad m(\xi_{n}):=\frac{1}{N}\sum_{1\leq i\leq N}\delta_{\xi^{i}_{n}}.
\end{equation}
This particle model (\ref{mean-field-intro}) is a genetic type particle model with a selection and a mutation transition dictated by the potential function $G_n$ and the Markov transition $M_{n}$. 

Loosely speaking, the approximation is constructed as follows: starting from a sample $\xi_0^{(N)}$ at $t=0$ of the initial distribution $\eta_0$ (so that $m(\xi_0)\simeq_{N\uparrow\infty}\eta_0$), and assuming
$m(\xi_n)\simeq_{N\uparrow\infty}\eta_n$, then the population at time $(n+1)$ 
is formed with  $N$ `almost' independent samples w.r.t. $\eta_{n+1}$ so that $m(\xi_{n+1})\simeq_{N\uparrow\infty}\eta_{n+1}$. The reader is refered to \cite{d-2004} for details.

We further assume that (\ref{FK-def-intro-ref}) is the historical version of the Feynman-Kac model $\eta^{\prime}_n$ associated with the parameters $(G_n^{\prime},M_n^{\prime})$,
 the $i$-th path space particle 
\begin{equation}\label{FK-path-particle}
\xi_n^i=\left(\xi^i_{0,n},\xi^i_{1,n},\ldots,\xi^i_{n,n}\right)\in S_n:=(S_0^{\prime}\times\ldots\times S^{\prime}_n)
\end{equation}
can be interpreted as
 the  line of ancestors $\xi^i_{p,n}$ of the $i$-th individual $\xi^i_{n,n}$ 
at time $n$, at every level $0\leq p\leq n$, with $1\leq i\leq N$. In addition the Markov chain $\xi^{\prime}_n:=\xi_{n,n}$ evolves as in (\ref{mean-field-intro})
by replacing $\Phi_n$ by the one step evolution semigroup $\Phi_n^{\prime}$ of the Feynman-Kac model $\eta^{\prime}_n$.

The path space model $\xi_n$ is called the genealogical tree model associated with the particle system $\xi_n$. 
We let $\XX_n\in S_n$ be a randomly chosen ancestral line with (conditional) distribution $m(\xi_n)$. We further assume that the integral operators $G^{\prime}_{k}(\textsl{x})M^{\prime}_{k+1}(\textsl{x}$ $,d\textsl{y})$ have a density $H^{\prime}_{k+1}(\textsl{x},\textsl{y})$ w.r.t. some reference measure $\nu^{\prime}_{k+1}$ on $S^{\prime}_{k+1}$. In this situation,  we have
\begin{equation}\label{df3}
\mbox{\rm Law}\left(\XX_n~|~(\xi^{\prime}_k)_{0\leq k\leq n}\right)=m(\xi^{\prime}_n)\HH_{n,m(\xi^{\prime})}
\end{equation}
with $m(\xi^{\prime})=(m(\xi^{\prime}_k))_{k\geq 0}$ and  the backward Markov transition $\HH_{n,m(\xi^{\prime})}$ from $S^{\prime}_n$ into $S_n$
defined by
$$
\HH_{n,m(\xi^{\prime})}(\textsl{x}_n,d(\textsl{y}_0,\ldots,\textsl{y}_n)):=\delta_{\textsl{x}_n}(d\textsl{y}_n)~\prod_{0\leq p<n}\frac{m(\xi^{\prime}_k)(d\textsl{y}_{k}) H^{\prime}_{k+1}(\textsl{y}_{k},\textsl{y}_{k+1})}{m(\xi^{\prime}_k)\left(H^{\prime}_{k+1}(\point,\textsl{y}_{k+1})\right)}.
$$
The proof of (\ref{df3}) can be found in~\cite{dkp-14} (see, for instance Proposition 4.6 and Theorem 4.7).

\subsection{An Upper and Lower Bound Bias Estimate}

\begin{hypA}\label{hyp:a1}
For any $n\geq 0$, there exist a $\rho_n<\infty$ s.t. for any $\textsl{x}_p,\textsl{y}_p\in S_p$, $0\leq p\leq q\leq n$
\begin{equation}\label{H}
Q_{p,q}(1)(\textsl{x}_p)\leq \rho_n~Q_{p,q}(1)(\textsl{y}_p).
\end{equation}
\end{hypA}

Define $\eta^{\textsl{x}}_{p,n}(\cdot)\propto Q_{p,n}(\textsl{x},\cdot)$  at time $p\leq n$, for $\textsl{x}\in S_p$.
%the solution of these equations starting at the Dirac measure $\delta_{\textsl{x}}$ 
\begin{hypA}\label{hyp:a2}
There exist a function $\alpha:\NN\rightarrow\mathbb{R}^+$ with
$
\sum_{n\geq 0}\alpha(n)<\infty,
$
such that for any $0\leq p\leq n$,
\begin{equation}\label{H-bis}
\sup_{\textsl{x},\textsl{y}}\left\Vert \eta^{\textsl{x}}_{p,n}-\eta^{\textsl{y}}_{p,n}\right\Vert_{\tiny tv}\leq \alpha(n-p).
%\quad \mbox{\rm for some function $\alpha$ s.t.}\quad\sum_{n\geq 0}\alpha(n)<\infty
\end{equation}
In addition, there exist a $0<\beta_1<\beta_2<\infty$, such that $0<\beta_1\leq \inf_{n,\textsl{x}} G_n(\textsl{x})< \sup_{n,\textsl{x}} G_n(\textsl{x})\leq \beta_2$.
\end{hypA}

\begin{rem}
In the context of (A\ref{hyp:a2}), we have
 $$
 (\ref{H-bis})~\Longrightarrow~
\rho:= \sup_{n\geq 0}\rho_n<\infty.
 $$
A proof of this result is provided in Chapter 12.2.1 in~\cite{d-2013}.
\end{rem}

By construction,
we have $Q_{p,q}(1)(\textsl{x}_p)=Q_{p,q}^{\prime}(1)(\textsl{x}_p^{\prime})$ for any $\textsl{x}_p=(\textsl{x}^{\prime}_k)_{0\leq k\leq p}\in S_p$; where $Q_{p,q}^{\prime}$ stands for the semigroup
defined as $Q_{p,q}$ by replacing $(S_n,G_n,M_n)$ by $(S^{\prime}_n,G^{\prime}_n,M^{\prime}_n)$. This implies that 
(A\ref{hyp:a1}-\ref{hyp:a2}) are met as soon as it is the case for the marginals.

%(\ref{H}) is met as soon as
%the $S^{\prime}_n$-marginal measures $\eta^{\prime,\textsl{x}}_{p,n}$ of $\eta^{\textsl{x}}_{p,n}$ satisfy the stability condition (\ref{H-bis}).

These conditions ensure that the mean field particle local sampling error do not propagate w.r.t.~the time horizon.
 They also guarantee that the bias and the variance of the occupation measures of path particles are $\mathcal{O}(n/N)$
  (see, for instance, Corollary 8.5.1 in~\cite{d-2004}, as well as Corollary 15.2.5 and Theorem 16.5.1 in \cite{d-2013}). In the context of particle Gibbs models
  with ancestral line sampling, under similar conditions these linear scaling were also obtained in~\cite{dkp-14,ldm-2014}.

\begin{rem}
 Replacing $G_n$ by  the normalized functions $G_n/\eta_n(G_n)$, there is no loss of generality 
to assume that the potential functions $G_n$ are chosen so that $\eta_n(G_n)=1$. In this situation,   we have
$$ \overline{Q}_{p,q}=Q_{p,q}\quad \mbox{and} \quad \eta_n=\gamma_n
$$ for any $0\leq p\leq q\leq n$ and
\begin{equation}\label{o1}
(\ref{H})\Longrightarrow \rho_n^{-1}\leq Q_{p,q}(1)(\textsl{x}_p)=\overline{Q}_{p,q}(1)(\textsl{x}_p)={Q_{p,q}(1)(\textsl{x}_p)}/{\eta_pQ_{p,q}(1)}\leq \rho_n
\end{equation}
\end{rem}

In this context, one of the main results of this article can be stated as follows.

\begin{theo}\label{T2}
Assume (A\ref{hyp:a1}). Then, for any $f:S_n\rightarrow\mathbb{R}^+$ s.t. $\eta_n(f)<\infty$ we have
{\small\begin{equation}\label{f15}
 \left(1-c_1(\rho_n) \frac{n}{N}\right)\eta_n(f) ~\leq~ 
\EE\left(m(\xi_{n})(f)\right)=\EE\left(m(\xi^{\prime}_n)\HH_{n,m(\xi^{\prime})}(f)\right) ~\leq~  \left(1+c_2(\rho_n)~\frac{n}{N}\right)\eta_n(f).
\end{equation}}
for any $N\geq c_2(\rho_n)n$ with
$$
c_1(\rho)=\rho^2-1\quad\mbox{\rm and} \quad c_2(\rho)= 2c_1(\rho) (2c_1(\rho)(c_1(\rho)+1)+1).
$$
In particular,
for any function $f\in \LL_1(\eta_n)$ and for any $N\geq c_2(\rho_n)n$
we have the bias estimates
\begin{equation}
\left\vert \EE\left(m(\xi_{n})(f)\right) -\eta_n(f)\right\vert=\left\vert \EE\left(m(\xi^{\prime}_n)\HH_{n,m(\xi^{\prime})}(f)\right) -\eta_n(f)\right\vert\leq c_2(\rho_n)~\frac{n}{N}~\eta_n(\vert f\vert).\label{eq:bias_est}
\end{equation}
\end{theo}

\begin{rem}
This Theorem  (\eqref{eq:bias_est}) establishes bias estimates for particle approximations and backward particle approximations for possibly unbounded functions on the path space $S_n$.
If one additionally assumes (A\ref{hyp:a2}), then these bounds are of $\mathcal{O}(n/N)$ (recall \eqref{H-bis}).
 This is an extremely important result that, to our knowledge, has not been previously established in the literature.
\end{rem}

\subsection{Sharp Propagation of Chaos Estimates}\label{propagation-chaos-sec}

This short section provides some comments on the proof of (\ref{f15}) and discusses some natural extensions of these estimates to tensor product empirical measures.
Roughly speaking, the proof of Theorem~\ref{T2} is based on the fact that
$$
\EE\left(
\frac{m(\xi_{p})(f)}{m(\xi_{p})(g)}
~\Big|~\xi_{p-1}\right)~
$$
$$
\simeq_{N\uparrow\infty}\frac{\EE(m(\xi_{p})(f)~|~\xi_{p-1})}{\EE(m(\xi_{p})(g)~|~\xi_{p-1})}~ \left(1+ \frac{c}{N}\right)=\frac{m(\xi_{p-1}) Q_p(f)}{m(\xi_{p-1}) Q_p(g)} ~\left(1+ \frac{c}{N}\right)
$$ 
for any non-negative functions $f,g$ on $S_p$ (s.t. $0<a\leq g\leq b$ for some $a,b>0$), and
for some finite constant $c<\infty$ whose values only depends on the parameters $(a,b)$. Iterating these inequalities backwards and recalling that $\eta_n(f)\propto \eta_0Q_{0,n}(f)$ we prove (\ref{f15}).

In the following, we require (A\ref{hyp:a1}-\ref{hyp:a2}) to hold.
One can replace the occupation measures $m(\xi_{p})$ by the $2$-tensor product measures with no replacement:
$$
m(\xi_p)^{\odot 2}=\frac{1}{N(N-1)}\sum_{1\leq i_1\not= i_2\leq N}\delta_{(\xi^{i_1}_p,\xi^{i_2}_p)}.
$$
Using similar techniques as described above, we can prove that
$$
\begin{array}{l}
\displaystyle\EE\left(
\frac{m(\xi_{p})^{\odot 2}(f)}{m(\xi_{p})^{\odot 2}(g)}
~\Big|~\xi_{p-1}\right)~\\
\\
\displaystyle\simeq_{N\uparrow\infty}~\frac{\EE(m(\xi_{p})^{\odot 2}(f)~|~\xi_{p-1})}{\EE(m(\xi_{p})^{\odot 2}(g)~|~\xi_{p-1})}~ \left(1+ \frac{c}{N}\right)
=\frac{m(\xi_{p-1})^{\otimes 2} Q_p^{\otimes 2}(f)}{m(\xi_{p-1})^{\otimes 2} Q_p^{\otimes 2}(g)}~ \left(1+ \frac{c}{N}\right)
\end{array}
$$
for any non-negative functions $f,g$ on $S^2_p$ (s.t. $0<a\leq g\leq b$ for some $a,b>0$), and 
for some finite constant $c<\infty$ whose values only depends on the parameter $\rho$. We recall  that 
$
m(\xi_p)^{\otimes 2}=m(\xi_p)^{\odot 2} \Ca
$,
with the coalescence operator
$
\Ca(f)(\textsl{x},\textsl{y}):=\left(1-\frac{1}{N}\right)~f(\textsl{x},\textsl{y})+\frac{1}{N}~  f(\textsl{x},\textsl{x})
$.
 In this situation, the backward recursion takes the form
$$
\EE\left(
\frac{m(\xi_{p})^{\odot 2}(f)}{m(\xi_{p})^{\odot 2}(g)}
~\Big|~\xi_{p-1}\right)~\simeq_{N\uparrow\infty}~\frac{m(\xi_{p-1})^{\odot 2}Q^{(2)}_p(f)}{m(\xi_{p-1})^{\odot 2} Q^{(2)}_p(g)} \left(1+ \frac{c}{N}\right)
$$
with $Q^{(2)}_p:=\Ca Q_p^{\otimes 2}$. Iterating these inequalities, we prove that
\begin{equation}\label{f17}
\left(1-a_1(\rho) \frac{n}{N}\right)\eta_n^{(2)}(f)\leq
\EE\left(m(\xi_{n})^{\odot 2}(f)\right)=\EE(f(\xi^1_n,\xi^2_n))\leq \left(1+a_2(\rho)~\frac{n}{N}\right)\eta^{(2)}_n(f).
\end{equation}
for some finite constants $a_1(\rho),a_2(\rho)<\infty$. In the above display, $\eta^{(2)}_n$ stands for the Feynman-Kac measures on $S_n^2$
associated with a reference Markov chain $X^{(2)}_n$ evolving on $S^2_n$ with Markov transitions $M_n^{(2)}:=M_n^{\otimes 2}\Ca$, initial distribution
$\eta^{(2)}_0:=\eta^{\otimes 2}_0\Ca$ and potential functions $G^{(2)}_n(\textsl{x},\textsl{y})=G_n(\textsl{x})\times G_n(\textsl{y})$. Following the arguments developed in~\cite{cdmg-11} we can check that
$$
\left(1-b_1(\rho) \frac{n}{N}\right)~\eta_n(\varphi)^2~\leq \eta_n^{(2)}(\varphi^{\otimes 2})=\EE(\varphi(\xi^1_n)\varphi(\xi^2_n))\leq \eta_n(\varphi)^2+b_2(\rho) \frac{n}{N} ~\eta_n(\varphi^2)
$$
for some finite constants $b_1(\rho),b_2(\rho)<\infty$. 

The extension of these estimates to $q$-tensor product measures can be developed
using the combinatorial and functional techniques developed in~\cite{dpr-2009,dkp-14}.
The technical details are available in Appendix \ref{tensor-section} and \ref{pce-section}.
 In these settings the $q$-tensor product empirical measure 
$m(\xi_p)^{\odot q}$ is defined as above by summing Dirac measures on $q$ particles  with different indices.
The coalescent operator $\Ca$ on $S^q_n$ connecting these measures to the conventional $q$-tensor product $m(\xi_p)^{\otimes q}=m(\xi_p)^{\odot q}\Ca$ is provided in
Proposition 8.6.1 in~\cite{d-2004} (see also Lemma~\ref{L3} in the Appendix). In this context, the estimates are described in terms of  the Feynman-Kac measures $\eta^{(q)}_n$  on $S_n^q$
associated with a reference Markov chain $X^{(q)}_n$ evolving on $S^q_n$ with Markov transitions $M_n^{(q)}:=M_n^{\otimes q}\Ca$, initial distribution
$\eta^{(q)}_0:=\eta^{\otimes q}_0\Ca$ and potential functions $G^{(q)}_n:=G_n^{\otimes q}$. 
In this notation, (\ref{f17}) remains valid if we replace 
$m(\xi_{n})^{\odot 2}$ and $\eta^{(2)}_n$ by $m(\xi_{n})^{\odot q}$ and $\eta^{(q)}_n$. In addition, one establishes sharp upper and lower bound
propagation of chaos estimates on path space:

\begin{theo}
Assume (A\ref{hyp:a1}-\ref{hyp:a2}). Then for 
any $1\leq q<N$, there exist $a(\rho)<\infty$
such that for any $n\geq 0$, any $f:S_n\rightarrow\mathbb{R}^+$ such that $\eta_n(f^q)<\infty$
 and $N>na(\rho)$, we have:
$$
\left(1-a(\rho)~\frac{n}{N}\right)~\eta_n(f)^q\leq 
\mathbb{E}\left[\prod_{1\leq i\leq q}f(\xi_n^i)\right]  \leq \eta_n(f)^q+a(\rho)\frac{n}{N}~\eta_n(f^q).
$$
\end{theo}

\section{Dual Processes and the Particle Gibbs Sampler}\label{sec:dual}

\subsection{Description of the Models}

We fix some random path $\textsl{z}:=(\textsl{z}_n)_{n\geq 0}\in \prod_{n\geq 0}S_n$. Notice that each coordinate
$\textsl{z}_n=\left(\textsl{z}_{k,n}\right)_{0\leq k\leq n}\in S_n=\prod_{0\leq k\leq n}S^{\prime}_k$ is itself a sequence of states $\textsl{z}_{k,n}$ in $S^{\prime}_k$, with $0\leq k\leq n$.

We consider the collection of Markov transitions 
$\Ma_{\textsl{z},n}$ from $S_{n-1}^N$ into $S_n^N$ defined by
$$
\Ma_{\textsl{z},n}(x_{n-1},dx_n)=\frac{1}{N}\left[\sum_{i=1}^{N}\Phi_{n}(m(x_{n-1})
)^{\otimes (i-1)}\otimes \delta_{\textsl{z}_{n}}\otimes\Phi_{n}(m(x_{n-1})
)^{\otimes (N-i)}\right](dx_n)
$$
and the distribution $\mu_{\textsl{z},0}$ on $S_0^N$ given by
$$
\mu_{\textsl{z},0}:=\frac{1}{N}\sum_{i=1}^{N}\left(\eta_0
^{\otimes (i-1)}\otimes \delta_{\textsl{z}_{0}}\otimes\eta_0^{\otimes (N-i)}\right).
$$

We let $\Xa_{\textsl{z},n}=(\Xa_{\textsl{z},n}^i)_{1\leq i\leq N}\in S_n^N$ be an $N$-particle model with initial distribution $\mu_{\textsl{z},0}$ and Markov transitions 
$\Ma_{\textsl{z},n}(x_{n-1},dx_n)$ from $S_{n-1}^N$ into $S_n^N$.
The process $\Xa_{X,n}$ associated with a frozen path $X=(X_n)_{n\geq 0}\in \prod_{n\geq 0}S_n$ is called the dual mean field model associated with the Feynman-Kac particle model $\xi_n$. 

We emphasize that each particle $\Xa_{\textsl{z},n}^i\in S_n\left(=\prod_{0\leq k\leq n}S^{\prime}_k\right)$ is a random path, with $1\leq i\leq N$. 
Given the system $\Xa_{\textsl{z},n-1}$ at rank $(n-1)$, we select randomly an index $I_n=i\in\{1,\ldots,N\}$ and we set $\Xa_{\textsl{z},n}^i=\textsl{z}_n$.
The other random variables $\Xa_{\textsl{z},n}^j$ with $j\in \{1,\ldots,N\}-\{i\}$ are $(N-1)$ conditionally independent random paths
with common law $\Phi_{n}\left(m(\Xa_{\textsl{z},n-1})\right)$.  They are defined as the genetic selection-mutation transition 
of genealogical tree models discussed in (\ref{mean-field-intro}) and (\ref{FK-path-particle}).

We set
$$
\zeta_n:=(\xi_k)_{0\leq k\leq n}\in \Sa_n:=\prod_{0\leq k\leq n}S_k^N\quad\mbox{\rm and}\quad \zeta^{\#}_{n}:=
(\Xa_{\textsl{X},k})_{0\leq k\leq n}\in \Sa_n.
$$
The main reason of the terminology "duality" comes from the fact that (see \cite{dkp-14}, e.g.~Theorem 4.3)
\begin{equation}\label{df1}
\EE\left(f_n\left(\XX_n,\zeta_n\right)~\Za_n(\zeta)\right)=\EE\left(f_n\left(X_n,\zeta^{\#}_{n}\right)~
Z_n(X)\right)\quad\mbox{\rm with}\quad \Za_n(\zeta):=\prod_{0\leq p<n}m(\xi_p)(G_p)
\end{equation}
for any bounded measurable function $f_n$ on $(S_n\times\Sa_n)$ (symmetric w.r.t. the coordinates in $S_k^N$ for each $0\leq k\leq n$). 
When (\ref{FK-def-intro-ref}) is the historical version of the Feynman-Kac model $\eta^{\prime}_n$ associated with the parameters $(G_n^{\prime},M_n^{\prime})$
we clearly have 
$$
Z_n(X)=\prod_{0\leq k<n}G^{\prime}_k(X^{\prime}_k):=Z_n^{\prime}(X^{\prime})\quad\mbox{\rm and}\quad
\Za_n(\xi)=\prod_{0\leq p<n}m(\xi_p^{\prime})(G^{\prime}_p):=\Za^{\prime}_n(\xi^{\prime}).
$$
In addition, formula (\ref{df1}) implies that (see again \cite{dkp-14}, e.g.~Theorem 4.3)
\begin{equation}\label{df2}
\EE\left(f_n\left(\XX_n,\zeta^{\prime}_n\right)~\Za_n^{\prime}(\zeta^{\prime})\right)=\EE\left(f_n\left(X_n,\zeta^{\prime \# }_{n}\right)~
Z_n^{\prime}(X^{\prime})\right)
\end{equation}
with
$$
\zeta^{\prime}_n:=(\xi_k^{\prime})_{0\leq k\leq n}\in  \Sa^{\prime}_n:=\prod_{0\leq p\leq n}S_p^{\prime N}
\quad\mbox{\rm and}\quad \zeta^{\prime \# }_{n}:=
(\Xa^{\prime}_{\textsl{X}^{\prime}\!,k})_{0\leq k\leq n}\in \Sa_n^{\prime}
$$
for any bounded measurable function $f_n$ on $(S_n\times\Sa_n^{\prime})$  (symmetric w.r.t.~the coordinates in $S_k^{\prime N}$ for each $0\leq k\leq n$). In the above display,
$\Xa^{\prime}_{\textsl{X}^{\prime},n}$ stands for the dual mean field model associated with the Feynman-Kac particle model $\xi^{\prime}_n$ and the frozen path $X^{\prime}=(X^{\prime}_n)_{n\geq 0}$.

%The proof of the duality formulae (\ref{df1}), (\ref{df2}) can be found in~\cite{dkp-14} (see for instance Theorem 4.3).

We let $\pi_n$ and $\pi_n^{\prime}$ be the probability measures on $(S_n\times\Sa_n)$ defined by
$$
\pi_n(f)\propto \EE\left(f_n\left(\XX_n,\zeta_n\right)~\Za_n(\zeta)\right)\quad\mbox{\rm and}\quad
\pi^{\flat}_n(f)\propto \EE\left(f_n\left(\XX_n^{\flat},\zeta^{\prime}_n\right)~\Za^{\prime}_n(\zeta^{\prime})\right)
$$
In the above display, $\XX_n^{\flat}$ stands for a backward random trajectory in $\in S_n$ with (conditional) distribution $m(\xi^{\prime}_n)\HH_{n,m(\xi^{\prime})}$.

The duality formulae (\ref{df1}) and (\ref{df2}) provide a natural way to design Gibbs-Glauber dynamics with reversible distribution $\pi_n$ and $\pi_n^{\flat}$ on the product space $\left(S_n\times\Sa_n \right)$ and $\left(S_n\times\Sa_n^{\prime}\right)$.
For instance,
the transition probabilities of the Gibbs-Glauber dynamics of the target multivariate distribution $\pi_n$
are described by the synthetic diagram
{\small $$
\left\{
\begin{array}{rcl}
\XX_n&=&\textsl{x}\\
\zeta_n&=&x
\end{array}
\right\}\longrightarrow \left\{
\begin{array}{rcl}
\overline{\XX}_n&=&\overline{\textsl{x}}\sim \left(\XX_n~|~\zeta_n=x\right)\\
\zeta_n&=&x
\end{array}
\right\}\longrightarrow
\left\{
\begin{array}{rcl}
\overline{\XX}_n&=&\overline{\textsl{x}}\\
\overline{\zeta}_n&=&\overline{x}\sim \left(\zeta_n~|~X_n=\overline{\textsl{x}}\right)
\end{array}
\right\}
$$}

In the above display, $\left(\XX_n~|~\zeta_n\right)$ and $\left(\zeta_n~|~X_n\right)$ 
is a shorthand notation for the $\pi_n$-conditional distributions of $\XX_n$ given $\zeta_n$, and $\zeta_n$ given $X_n$.
Notice that the first transition of the Gibbs-Glauber dynamics reduces to the uniform sampling of an ancestral line. In addition, by  (\ref{df1}), 
the second transition amounts of sampling a genetic particle model with a frozen ancestral line. One may interpret the above as a version of the particle
Gibbs sampler of \cite{adh-2010}.

The Markov transitions of the first coordinate
of the Gibbs-Glauber dynamics with ancestral path and backward sampling are defined for any $f\in\Ba(S_n)$ and $\textsl{z}\in S_n$ by
the formulae
$$
\KK_n(f)(\textsl{z}):=\EE\left(m(\Xa_{\textsl{z},n})(f)\right)\quad\mbox{\rm and}\quad
\KK^{\flat}_n(f)(\textsl{z}):=\EE\left(m(\Xa_{\textsl{z},n}^{\prime})\HH_{n,m(\Xa_{\textsl{z}}^{\prime})}(f)\right).
$$

The dual process $\Xa_{\textsl{z},n}$ can be interpreted as a mean field particle model with a given frozen path which also has a negligible
impact of $\mathcal{O}(1/N)$ upon the auxiliary particle system. This indicates  that both the particle density profiles $m(\Xa_{\textsl{z},n})(f)$, and the Feynman-Kac
measures $\eta_{\textsl{z},n}$ with a frozen path converge to $\eta_n$ as $N$ tends to $\infty$ (see, for instance, Proposition 4.12 in~\cite{dkp-14}). In addition, as shown in Proposition 4.11~\cite{dkp-14}, 
the first moment of the unnormalized density profiles of the frozen process coincides with the Feynman-Kac model with a frozen path; that is, we have that
$$
\gamma_{\textsl{z},n}(f):=
\EE\left(m(\Xa_{\textsl{z},n})(f)~\prod_{0\leq k<n}m(\Xa_{\textsl{z},k})(G_k)\right)
$$
where $(\gamma_{\textsl{z},n},\eta_{\textsl{z},n})$ stand for  the Feynman-Kac measures defined as  (\ref{FK-def-intro-ref}) by replacing 
$\eta_0=\mbox{\rm Law}(X_0)$ and $M_n$  by  the initial distribution $\eta_{\textsl{z},0}=\mu_{\textsl{z},0}$
 and
the Markov transitions $M_{\textsl{z},n}$ given by
\begin{eqnarray*}
M_{\textsl{z},n}(\textsl{x}_{n-1},d\textsl{x}_n)&=&\frac{1}{N}~\delta_{\textsl{z}_n}(d\textsl{x}_n)+\left(1-\frac{1}{N}\right)~M_{n}(\textsl{x}_{n-1},d\textsl{x}_n).
\end{eqnarray*}

We also consider the Markov transition $\FF_n$ from $S_n$ into itself defined by
the Feynman-Kac model with a frozen path; that is,
$$
\FF_n(f)(\textsl{z}):=\eta_{\textsl{z},n}(f):=\gamma_{\textsl{z},n}(f)/\gamma_{\textsl{z},n}(1).
$$

\subsection{A Comparison Theorem and Contraction Inequalities}

A first order expansion of  $(\gamma_{\textsl{z},n},\eta_{\textsl{z},n})$ around their limiting measures  $\gamma_{n}=\eta_{n}$ as $N\rightarrow\infty$
is provided below (recall that the potential functions $G_n$ are centered in the sense that $\eta_n(G_n)=1$).
\begin{prop}\label{P2}
%We assume that the regularity condition (\ref{H}) is satisfied for some finite $\rho<\infty$.
Assume (A\ref{hyp:a1}-\ref{hyp:a2}). Then, for any $f\in\Ba(S_n)$, $f:S_n\rightarrow[0,1]$ and  for any $N\geq 2(1+2\rho^2)n$ we have
\begin{equation}\label{f11}
 \left(1-2~\frac{n}{N}\right)~\eta_{n}(f)\leq \gamma_{\textsl{z},n}(f)\leq \eta_{n}(f)+4\rho^2~\frac{n}{N}
\end{equation}
and 
\begin{equation}\label{f6}
\left(1-2(1+2\rho^2)~\frac{n}{N}\right)~\eta_{n}(f)\leq \eta_{\textsl{z},n}(f)~\leq \eta_{n}(f)+4(1+2\rho^2)~\frac{n}{N}.
\end{equation} 
\end{prop}

The formulas (\ref{f6}) and (\ref{f11}) are direct consequences of the semigroup estimates stated in Proposition~\ref{P1} in Section \ref{sec:semi}. 

By construction, the Feynman-Kac distribution $\eta_n$ on path space  is the unique equilibrium measure of $\KK_n$ and $\KK^{\flat}_n$, for any choice of the population size $N$ and for any time horizon $n$. In addition, we have
\begin{equation}\label{f14}
\sup_{\textsl{z}\in S_n}{\left(\left\Vert \KK_n^p(\textsl{z},\point)-\eta_n\right\Vert_{\tiny tv}\vee \left\Vert \KK_n^{\flat, p}(\textsl{z},\point)-\eta_n\right\Vert_{\tiny tv}\right)}\leq (1-\epsilon_n(N))^p~
\end{equation}
for some finite constant $\epsilon_n^N\in ]0,1]$ whose values only depends on the parameters $(n,N)$.

For instance, let us assume that the potential functions are chosen s.t.
 $0<\epsilon\leq G_k\leq 1$, for some $\epsilon\in ]0,1]$, with $k\geq 0$. In this situation, for any $\textsl{x}_n=(\textsl{x}_{n-1},\textsl{x}^{\prime}_n)\in
S_n=\left(S_{n-1}\times S^{\prime}_n\right)$ and any non-negative function $f$ on $S_n$ we have
the rather crude minorization estimate
$$
\KK_n^{\flat}(f)(\textsl{x}_n)\geq  \epsilon( 1-1/N)~\KK_{n-1}^{\flat}(\overline{Q}_n(f))(\textsl{x}_{n-1})\Rightarrow
\KK_n^{\flat}(f)(\textsl{x}_n)\geq (\epsilon( 1-1/N))^n~\eta_n(f)
$$
as soon as $N\geq 2$. In much the same way, we prove that
\begin{eqnarray*}
\KK_n(f)(\textsl{x}_n)
&=&\epsilon( 1-1/N)~\KK_{n-1}(\overline{Q}_n(f))(\textsl{x}_{n-1})\Rightarrow
\KK_n(f)(\textsl{x}_n)\geq (\epsilon( 1-1/N))^n~\eta_n(f).
\end{eqnarray*}
These couple of estimates imply (\ref{f14}) with $$
\epsilon_n(N)\geq (\epsilon( 1-1/N))^n\geq \epsilon^n~(1-n/N)
$$
The r.h.s.~estimate comes from the fact that $x^q=1-(1-x^q)$ and $(1-x^q)=(1-x)\sum_{0\leq k<q}x^k\leq q(1-x)$ for any $0\leq x<1$. 
For a more detailed discussion on these rather crude estimates, we refer the reader to Section 4.3, Proposition 4.10 in~\cite{dkp-14}.

The following theorem, which is proved using the same line of arguments as Theorem~\ref{T2}, discusses some rather mild conditions under which $\epsilon_n(N)\geq 1-cn/N$ for any $(n,N)$ such that $N\geq cn$,
for some finite constant $c<\infty$.

\begin{theo}\label{T1}
%We assume that the regularity condition (\ref{H}) is satisfied for some finite $\rho<\infty$.
Assume (A\ref{hyp:a1}-\ref{hyp:a2}). Then, there exists some finite constant $c(\rho)<\infty$ such that
for any $f\in\Ba(S_n)$, $f:S_n\rightarrow\mathbb{R}^+$ and any $N\geq c(\rho) n$ 
we have
\begin{equation}\label{f5}
\left(1-c(\rho)~\frac{n}{N}\right)~\FF_n(f)(\textsl{z}):~\leq \KK_n(f)(\textsl{z})~\leq 
\left(1+c(\rho) ~\frac{n}{N}\right)~\FF_n(f)(\textsl{z})
\end{equation}
and
\begin{equation}\label{f12}
\KK^{\flat}_n(f)(\textsl{z})~\geq  \left(1-c(\rho)~\frac{n}{N}\right)~\eta_{n}(f)/\gamma_{\textsl{z},n}(1)
\end{equation}
 
\end{theo}
The estimates (\ref{f5}) and (\ref{f12}) are proved in section \ref{sec:prf_theo1}.

Combining Theorem~\ref{T1} with Proposition~\ref{P2}, and using the fact that 
\begin{equation}\label{f24}
\left[(1-x)/(1+y)\right]\vee \left[(1-x)(1-y)\right]\geq 1-(x+y)
\end{equation} for any $0<x,y\leq 1$, we conclude that
\begin{equation}
N\geq c(\rho)n\Longrightarrow \KK_n(f)(\textsl{z})\wedge \KK^{\flat}_n(f)(\textsl{z})~\geq \left(1-\frac{c(\rho)n}{N}\right)~\eta_n(f)
\label{eq:minor}
\end{equation}
for any non-negative function $f\in\Ba(S_n)$, and some  finite constant $c(\rho)<\infty$.

As a direct consequence of the estimates (\ref{f6}) and (\ref{f5}), for any function $f\in\Ba(S_n)$ s.t.
 $\mbox{\rm osc}(f)\leq 1$ we have 
\begin{equation}\label{f16}
\begin{array}[b]{l}
\displaystyle\mbox{\rm osc}\left(\FF_n(f)\right)\leq c_1(\rho)~\frac{n}{N}\quad
\mbox{\rm and}\quad
\mbox{\rm osc}\left(\KK_n(f)\right)\leq \mbox{\rm osc}\left(\FF_n(f)\right)+c_2(\rho)~\frac{n}{N}\\
\\
\displaystyle\Longrightarrow \mbox{\rm osc}\left(\KK_n(f)\right)\leq [c_1(\rho)+c_2(\rho)]~\frac{n}{N}
\end{array}
\end{equation}
for some finite constants $c_i(\rho)<\infty$, $i=1,2$, as soon as $N\geq c(\rho):=c_1(\rho)+c_2(\rho)$. In terms of 
 Dobrushin contraction coefficient (\ref{f27}) we have proved that
$$
\beta(\FF_n)\leq c_1(\rho)~\frac{n}{N}
\quad\mbox{\rm and}\quad
\beta(\KK_n)\leq \beta(\FF_n)+c_2(\rho)~\frac{n}{N}\Rightarrow \beta(\KK_n)\leq c(\rho)~\frac{n}{N}
$$

The ideas in this section are inspired by the work of \cite{alv-2013,dkp-14,ldm-2014}.
The minorization conditions (\ref{eq:minor}) are in a similar vein to those in \cite{alv-2013,ldm-2014}, in that we only establish uniform geometric convergence to
$\eta_n$. We also provide an explicit proof for backward sampling, although \cite{ldm-2014} state that their results extend to such a case. The rates in 
\eqref{eq:minor} are also established in \cite{alv-2013,ldm-2014} albeit using a different technical approach. Note, however, that they do not have an expansion of the kernel
as in \eqref{f5}.
The oscillation estimates (\ref{f16}) are also in similar vein to those in \cite{dkp-14} using non-asymptotic Taylor expansions of the transition $\KK_n$ around the equilibrium measure $\eta_n$. Our approach is based on
two elementary first order estimates (\ref{f6}) and (\ref{f5}) relating particle Gibbs transitions to Feynman-Kac models on path space with a reference Markov chain
that may jump from time to time to the frozen ancestral line.

We note that the result is extended to the scenario
of the inclusion of a static parameter $\theta$, assuming that the Gibbs update on $\theta$ possess a one-step minorizarion condition and (A\ref{hyp:a1}-\ref{hyp:a2}) holds
uniformly in $\theta$: the kernel which first updates $\theta$ from its full conditional and then applies the kernel(s) analyzed above then possess the above properties.

\subsection{Propagation of Chaos Estimates}

The first order propagation of chaos estimates can be derived extending the propagation of chaos analysis developed in Section~\ref{propagation-chaos-sec}
to empirical measures of the dual particle model with a frozen path. To precisely describe these extensions, we let 
$Q_n^{(q)}:=\Ca Q_n^{\otimes q}$ be the Feynman-Kac semigroup  associated with the Feynman-Kac measures $\eta_n^{(q)}$ discussed in the end of 
Section~\ref{propagation-chaos-sec} (a more thorough discussion on these Feynman-Kac models is provided in Appendix~\ref{pce-section}). 
We also consider the Markov operator $\Da_{\textsl{z}_n}$ from $S_n^q$ into itself defined for any function $f\in \Ba(S^q_n)$ by
$$
\Da_{\textsl{z}_n}(f)(\textsl{x}^1,\ldots,\textsl{x}^q):=\frac{1}{q}\sum_{1\leq i\leq q}f(\textsl{x}^1,\ldots,\textsl{x}^{i-1},\textsl{z}_n,\textsl{x}^{i+1},\ldots,\textsl{x}^q).
$$

In the above notation, the central idea it to observe that
\begin{equation}\label{f28}
\EE\left(m(\Xa_{\textsl{z},n})^{\odot q}(f)~|~\Xa_{\textsl{z},n-1}~\right)=\frac{ m(\Xa_{\textsl{z},n-1})^{\odot q} Q_{\textsl{z},n}^{\{q\}}(f)}{m(\Xa_{\textsl{z},n-1})^{\odot q} Q_{\textsl{z},n}^{\{q\}}(1)}
\quad\mbox{\rm with}\quad
Q_{\textsl{z},n}^{\{q\}}:=Q_{n}^{(q)}\Da_{\textsl{z}_n}.
\end{equation}
 To check this claim, we recall that
 $$
 \EE\left(m(\Xa_{\textsl{z},n})^{\odot q}(f)~|~\Xa_{\textsl{z},n-1}~\right)=\frac{1}{(N)_q}\sum_{a\in I^N_q} \EE\left(
 f(\Xa_{\textsl{z},n}^{a(1)},\ldots,\Xa_{\textsl{z},n}^{a(q)})~|~\Xa_{\textsl{z},n-1}~\right)
 $$
where $I^N_q$ stands for the set of all the $(N)_q=N!/(N-q)!$ one to one mappings from $\{1,\ldots,q\}$ into
 $\{1,\ldots,N\}$. By construction, we have
 $$
 \begin{array}{l}
 \EE\left(
 f(\Xa_{\textsl{z},n}^{a(1)},\ldots,\Xa_{\textsl{z},n}^{a(q)})~|~\Xa_{\textsl{z},n-1}~\right)\\
 \\
\displaystyle =\frac{1}{N}\sum_{1\leq i\leq N}\int \left\{\prod_{1\leq j\not=i\leq N}\Phi_{n}(m(\Xa_{\textsl{z},n-1}))(d\textsl{x}_n^j)\right\}~\delta_{\textsl{z}_n}(d\textsl{x}^i_n)~
 f(\textsl{x}_n^{a(1)},\ldots,\textsl{x}_n^{a(q)})\\
 \\
\displaystyle =\frac{q}{N}~\Phi_{n}(m(\Xa_{\textsl{z},n-1}))^{\otimes q}(\Da_{\textsl{z}_n}(f))+\left(1-\frac{q}{N}\right)~\Phi_{n}(m(\Xa_{\textsl{z},n-1}))^{\otimes q}(f).
\end{array}
$$

We let $Q^{\{q\}}_{\textsl{z},p,n}=Q^{\{q\}}_{\textsl{z},p+1}Q^{\{q\}}_{\textsl{z},p+1,n}$ be the semigroup associated with
the integral operator $Q^{\{q\}}_{\textsl{z},n}$, and we consider the corresponding Feynman-Kac measures
$$
\eta^{\{q\}}_{\textsl{z},n}(f):=\gamma^{\{q\}}_{\textsl{z},n}(f)/\gamma^{\{q\}}_{\textsl{z},n}(1)\quad\mbox{\rm with}\quad \gamma^{\{q\}}_{\textsl{z},n}:=\eta^{\{q\}}_{\textsl{z},0}Q^{\{q\}}_{\textsl{z},0,n} \quad \mbox{\rm and}\quad \eta^{\{q\}}_{\textsl{z},0}:=\eta_0^{\otimes q}\Da_{\textsl{z}_0}
$$
Arguing as in the proof of Proposition~\ref{P1}, under the assumptions (A\ref{hyp:a1}-\ref{hyp:a2}) we check that
$$
\left(1-a_q(\rho)\frac{n}{N}\right)
\eta^{(q)}_{n}(f)
\leq 
\eta^{\{q\}}_{\textsl{z},n}(f)\leq \eta^{(q)}_{n}(f)+ a_q(\rho)\frac{n}{N}
$$
for any $[0,1]$-valued function $f\in\Ba(S^q_n)$, some finite constants $a_q(\rho)$, and as soon as $N\geq n a_q(\rho)$. 
In the above display, $\eta^{(q)}_{n}$ are the Feynman-Kac measures discussed in the end of Section~\ref{propagation-chaos-sec} (see also Appendix~\ref{pce-section}).
In addition, extending the propagation of chaos analysis for tensor product measures $m(\xi_n)^{\odot q}$ to the empirical
measures $m(\Xa_{\textsl{z},n})^{\odot q}$ of the dual particle model, we prove the estimates
$$
\left(1-b(\rho)\frac{n}{N}\right)
 \eta^{\{q\}}_{\textsl{z},n}(f)
\leq \EE\left(m(\Xa_{\textsl{z},n})^{\odot q}(f)\right)\leq  \eta^{\{q\}}_{\textsl{z},n}(f)\left(1+b(\rho)\frac{n}{N}\right)
$$
for some finite constant $b(\rho)$, and as soon as $N\geq n b(\rho)$. 

\section{Semigroup Estimates}\label{sec:semi}

In this section we consider some original semigroup analysis, which leads to the proofs of Theorems \ref{T2} and \ref{T1}.

\subsection{Semigroup Analysis}

In this section we fix the time horizon $n$ and the frozen trajectory $\textsl{z}\in S_n$.
To simplify the presentation when there is no possible confusion we suppress the index $(\point)_{\textsl{z}}$.
For instance, we write $\Xa_{n}$ and $\Xa_n^{\prime}$  instead of $\Xa_{\textsl{z},n}$ and $\Xa_{\textsl{z},n}^{\prime}$. We also assume (A\ref{hyp:a1}-\ref{hyp:a2}). When the stability property (\ref{H-bis}) is not met, all the results developed in this section remain valid by replacing $\rho$ by $\rho_n$.

Given $N\geq 1$ and $\epsilon=\left(\epsilon_n\right)_{n\geq 1}\in \{0,1\}^{\NN}$, we consider the integral operators
\begin{eqnarray*}
\Qa_n&=&Q_n\Ia_n \quad\mbox{\rm with}\quad \Ia_n(\textsl{x}_n,d\textsl{y}_n):=\frac{1}{N}~\delta_{\textsl{z}_n}(d\textsl{y}_n)+\left(1-\frac{1}{N}\right)~\delta_{\textsl{x}_n}(d\textsl{y}_n)\\
\Qa_n^{(\epsilon)}&=&Q_n\Ia_n^{(\epsilon)} \quad\mbox{\rm with}\quad \Ia_n^{(\epsilon)}(\textsl{x}_n,d\textsl{y}_n):=\epsilon_n~\delta_{\textsl{z}_n}(d\textsl{y}_n)+\left(1-\epsilon_n\right)~\delta_{\textsl{x}_n}(d\textsl{y}_n)
\end{eqnarray*}
and the measure $\mu_0:=\eta_0\Ia_0$.
We denote by $\Qa_{p,n}$ and $\Qa_{p,n}^{(\epsilon)}$  the corresponding semigroup defined by the backward recursion 
$\Qa_{p-1,n}=\Qa_{p}\Qa_{p,n}$, and $\Qa_{p-1,n}^{(\epsilon)}=\Qa_{p}^{(\epsilon)}\Qa_{p,n}^{(\epsilon)}$  for $1\leq p\leq n$; 
with the convention $\Qa_{n,n}=\Qa^{(\epsilon)}_{n,n}=Id$, the identity integral operator on $\Ba(S_n)$ for $p=n$.

By construction, for any $p,q\geq 0$ we have
$$
\Qa_{p,p+q}=\sum_{0\leq k\leq q}\frac{1}{N^k}\left(1-\frac{1}{N}\right)^{q-k}\sum_{\epsilon_1+\ldots+\epsilon_q=k} \Qa_{p,p+q}^{(\epsilon)}.
$$
The r.h.s. summation in the above display is taken over all sequences $(\epsilon_i)_{1\leq i\leq q}\in\{0,1\}^q$ s.t. $\sum_{1\leq i\leq q}\epsilon_i=k$.
For instance, when $\epsilon_{p_i}=1$ with $1\leq i\leq k$, for some sequence of time steps $p+1\leq r_1<\ldots<r_k\leq p+q$, for any $f\in\Ba(S_{p+q})$ and $\textsl{x}_p\in S_p$ we have
$$
\Qa_{p,p+q}^{(\epsilon)}(f)(\textsl{x}_p)=
Q_{p,r_1-1}(1)(\textsl{x}_p)  \left\{\prod_{1\leq l< k}Q_{r_l,r_{l+1}-1}(1)(\textsl{z}_{r_l})\right\}
Q_{r_k,p+q}(f)(\textsl{z}_{r_k}).
$$

\begin{prop}\label{P1}
For any $p,q\geq 0$, $f\in \Ba(S_p)$ s.t. $0\leq f\leq 1$, and any $N\geq 3q\rho$ we have
\begin{equation}\label{i1}
 \left(1-\frac{q}{N}\right)~Q_{p,p+q}(f)\leq  \Qa_{p,p+q}(f)\leq Q_{p,p+q}(f)+\frac{2q\rho^2}{N}.
\end{equation}
In particular, for any  $0\leq p\leq n$ and any $N\geq 3n\rho$ we have the estimates
\begin{equation}\label{i2}
\frac{2}{3}~\rho^{-1}\leq \left(1-\frac{n}{N}\right)\rho^{-1}~\leq  \Qa_{p,n}(f)\leq  \left(1+\frac{2n\rho}{N}\right)~\rho\leq \frac{5}{3}~\rho.
\end{equation}
In addition,  for any $N\geq (1+2\rho^2)(n+1)$ we have
\begin{equation}\label{i3}
 \left(1-\frac{n+1}{N}\right)~\eta_{n}(f)\leq  \mu_0\Qa_{0,n}(f)\leq \eta_{n}(f)+2\rho^2~\frac{(n+1)}{N}
\end{equation}
and
\begin{equation}\label{i4}
 \left(1-(1+2\rho^2)~\frac{n+1}{N}\right)~\eta_{n}(f)\leq  \frac{\mu_0\Qa_{0,n}(f)}{ \mu_0\Qa_{0,n}(1)}\leq \eta_{n}(f)+2(1+2\rho^2)~\frac{2(n+1)}{N}.
\end{equation}
\end{prop}
\proof
Under our assumptions, we have 
$$
\sum_{1\leq i\leq q}\epsilon_i=k~\Longrightarrow~
\left\Vert
\Qa_{p,p+q}^{(\epsilon)}(1)
\right\Vert\leq \rho^{k+1}
$$
from which we prove that
\begin{eqnarray*}
 \Qa_{p,p+q}(f)(\textsl{x}_p)&\leq& Q_{p,p+q}(f)(\textsl{x}_p)+\rho~\left[
\left(\left(1-\frac{1}{N}\right)+\frac{\rho}{N}\right)^{q}-\left(1-\frac{1}{N}\right)^q
\right]\\
&\leq & Q_{p,p+q}(f)(\textsl{x}_p)+2q\rho^2/N
\end{eqnarray*}
for any function $f$ s.t. $0\leq f\leq 1$, and as soon as $N\geq 3q\rho$. The r.h.s.~estimate comes from the fact that
$$
\left(1+\frac{x}{N}\right)^{n}-\left(1-\frac{1}{N}\right)^{n}\leq 2(x+1)~\frac{n}{N}
$$
for any $x,n\geq 0$ and any $N\geq 3nx$. In the last estimate, we have used the 
standard decomposition 
$$
a^n-b^n=(a-b)\sum_{0\leq k<n}a^kb^{(n-1)-k}\leq n (a-b) a^n\leq 2n(a-b)
$$ 
as soon as $1\leq a^n\leq 2$ and $b\leq 1$. On the other hand, we have
$$
 \Qa_{p,p+q}(f)(\textsl{x}_p)\geq Q_{p,p+q}(f)(\textsl{x}_p)~\left(1-\frac{1}{N}\right)^q\geq Q_{p,p+q}(f)(\textsl{x}_p)~\left(1-\frac{q}{N}\right).
$$
This ends the proof of (\ref{i1}). 
The estimate (\ref{i2}) is a direct consequence of (\ref{i1}) with (\ref{o1}). To check (\ref{i3}) we observe that
$$
\left(1-\frac{(n+1)}{N}\right) \eta_n(f)\leq \left(1-\frac{1}{N}\right)^{n+1}~\eta_0Q_{0,n}(f)\leq  \mu_0\Qa_{0,n}(f).
 $$
In much the same way, we have
$$
 \mu_0\Qa_{0,n}(f)\leq \mu_0Q_{0,n}(f)+\frac{2n\rho^2}{N}\leq \eta_0Q_{0,n}(f)+\frac{2n\rho^2+\rho}{N}\leq \eta_n(f)+\frac{2(n+1)\rho^2}{N}.
$$
Now we come to the proof of the estimate (\ref{i4}). By (\ref{i3}) we have
$$
\frac{\mu_0\Qa_{0,n}(f)}{ \mu_0\Qa_{0,n}(1)}\geq ~\eta_{n}(f)~\frac{ \left(1-\frac{n+1}{N}\right)}{\left(1+2\rho^2~\frac{(n+1)}{N}\right)}\geq  
\eta_{n}(f)~\left(1-(2\rho^2+1)~\frac{n+1}{N}\right).
$$
 On the other hand, combining  (\ref{i3}) 
with the fact that
$$
\mu_0\Qa_{0,n}(1)\geq \left(1-\frac{1}{N}\right)^{n+1}\eta_0Q_{0,n}(1)=\left(1-\frac{1}{N}\right)^{n+1}\geq 1-\frac{n+1}{N}\geq \frac{1}{2}
$$
for $N\geq 2(n+1)$, we also have
$$
\frac{\mu_0\Qa_{0,n}(f)}{ \mu_0\Qa_{0,n}(1)}\leq \frac{1}{1-\frac{n+1}{N}}~\eta_{n}(f)+4\rho^2~\frac{(n+1)}{N}\leq \left(1+2~\frac{n+1}{N}\right)\eta_{n}(f)+4\rho^2~\frac{(n+1)}{N}.
$$
In the r.h.s.~estimate we have used the fact that $1/(1-x)\leq 1+2x$, for any $0<x\leq 1/2$. 
This ends the proof of the proposition.
\cqfd

\subsection{Technical Result}

The proof of Theorem~\ref{T2}, and the proof of the estimates (\ref{f5}) and (\ref{f12}) stated in Theorem~\ref{T1} are based on the following technical lemma of separate interest. 

\begin{lem}\label{L2}
Let $X=(X^1,\ldots,X^N)$ be a collection of $N$ conditionally independent random variables w.r.t.~some $\sigma$-field $\Fa$,
 taking values in some measurable state space $(E,\Ea)$, and such that
$$
\forall f\in\Ba(E)\qquad
\EE\left(m(X)(f)|\Fa\right)=\EE\left(m(X)(f)\right).
$$
For any $f,g\in \Ba(E)$, s.t. $0<a\leq g\leq b<\infty$, $f:E\rightarrow\mathbb{R}^+$, for some finite constant $(a,b)$
we have
\begin{equation}\label{f3}
\left(1-\frac{c_1}{N}\right)
~\frac{\EE\left(m(X)(f)\right)}{\EE\left(m(X)(g)\right)}
~\leq \EE\left(
\frac{m(X)(f)}{m(X)(g)}\right)\leq 
\frac{\EE\left(m(X)(f)\right)}{\EE\left(m(X)(g)\right)}~\left(1+\frac{c_2}{N}\right)
\end{equation}
as soon as $N>2(c_1+1)$, with
$$
c_1=\frac{b}{a}-1\quad\mbox{and}\quad c_2:= 2 c_1~(2c_1(c_1+1)+1)
$$
\end{lem}
Lemma~\ref{L2} is a direct consequence of the Proposition~\ref{p1} presented in the Appendix.

\subsection{Proof of Theorem~\ref{T2}}

{\bf Proof of Theorem~\ref{T2}}

We let $\Ha_p=\sigma(\xi_q,~q\leq p)$ be the natural $\sigma$-algebra filtration associated with the mean field particle model $\xi_{q}$, with $0\leq q\leq p$.
Notice that
$$
\EE\left(m(\xi_{p+1})(f)~|~\Ha_{p}\right)=\frac{m(\xi_{p})(Q_{p+1}(f))}{m(\xi_{p})(Q_{p+1}(1))}.
$$
Therefore, recalling that $Q_{p,n}=Q_{p+1}Q_{p+1,n}$, it is readily checked that
$$
\frac{\EE\left[ m(\xi_{p+1})(Q_{p+1,n}(f))|\Ha_p
\right]
}{\EE\left[ m(\xi_{p+1})(Q_{p+1,n}(1))|\Ha_p
\right]
}=\frac{m(\xi_{p})(Q_{p,n}(f))}{m(\xi_{p})(Q_{p,n}(1))}.
$$
Applying Lemma~\ref{L2} to the collection  of $N$ conditionally independent random variables $(\xi_{p+1}^{i})_{1\leq i\leq N}$ w.r.t. $\Ha_p$
and using the estimates (\ref{o1}) we find that
$$
\left(1-\frac{c_1(\rho)}{N}\right)\frac{m(\xi_{p})(Q_{p,n}(f))}{m(\xi_{p})(Q_{p,n}(1))}~\leq \EE\left[
\frac{m(\xi_{p+1})(Q_{p+1,n}(f))}{m(\xi_{p+1})(Q_{p+1,n}(1))}|\Ha_p
\right]~$$
and
$$
 \EE\left[
\frac{m(\xi_{p+1})(Q_{p+1,n}(f))}{m(\xi_{p+1})(Q_{p+1,n}(1))}\Big|\Ha_p
\right]~\leq 
\left(1+\frac{c_2(\rho)}{N}\right)\frac{m(\xi_{p})(Q_{p,n}(f))}{m(\xi_{p})(Q_{p,n}(1))}
$$
for any non-negative functions $f\in \Ba(S_n)$,
with $$c_1(\rho)=\rho^2-1\quad\mbox{\rm and} \quad c_2(\rho)= 2c_1(\rho) (2c_1(\rho)(c_1(\rho)+1)+1)$$ Taking expectations and iterating this process we prove that
\begin{eqnarray*}
\EE\left(m(\xi_{n})(f)\right)=\EE\left(\frac{m(\xi_{n-1})(Q_{n}(f))}{m(\xi_{n-1})(Q_{n}(1))}\right) & \leq &
\left(1+\frac{c_2(\rho)}{N}\right)^{n}\frac{\eta_0(Q_{0,n}(f))}{\eta_0(Q_{0,n}(1))}\\ & = &\left(1+\frac{c_2(\rho)}{N}\right)^{n}\eta_n(f).
\end{eqnarray*}
In much the same way, we have the lower bound estimate
$$\left(1-\frac{c_1(\rho)}{N}\right)^{n}\eta_n(f)\leq
\EE\left(m(\xi_{n})(f)\right).
$$

We check (\ref{f15}) recalling that $$x^n=1-(1-x^n)=1-(1-x)~\sum_{0\leq k<n}x^k\geq 1-n(1-x)$$ for any $0\leq x\leq 1$, and
$$
N\geq eny \Rightarrow\begin{array}[t]{rcl}
(1+y/N)^n-1=e^{n\log (1+y/N)}-1&\leq& e^{ny/N}-e^0\\
&\leq& ~e^{ny/N}~{ny}/{N}\leq {eny}/{N}\left(\leq 1\right).
\end{array}$$
This ends the proof of the estimates  (\ref{f15}). Notice that
\begin{eqnarray*}
(\ref{df3})\Longrightarrow \EE\left(m(\xi_{n})(f)\right)=\EE\left(f(\XX_n)\right) & = &\EE\left(\EE\left(f(\XX_n)\left\vert (\xi^{\prime}_k)_{0\leq k\leq n}\right.\right)\right)\\&=&\EE\left(m(\xi^{\prime}_n)\HH_{n,m(\xi^{\prime})}(f)\right).
\end{eqnarray*}
The proof of Theorem~\ref{T2} is now completed.

\subsection{Proof of Theorem \ref{T1}} \label{sec:prf_theo1}

{\bf Proof of (\ref{f5}) :}\label{pf5}

We let $\Ga_n=\sigma(\Xa_p,~p\leq n)$ be the natural $\sigma$-algebra filtration associated with the dual particle model $\Xa_{n}$.
Given $\Ga_{n}$, $\Xa_{n+1}$ are conditionally independent random variables 
w.r.t. $\Ga_{n}\wedge\Fa_n$, where $\Fa_n$ stands for the $\sigma$-algebra generated by an uniform random variable $I_n$ on $[N]$ (independent of $\Ga_{n}$).
More precisely, the conditional distribution of $\Xa_{n+1}$ w.r.t. $\Fa_{n}$ is 
given by
$$
\mbox{\rm Law}\left(\Xa_{n+1}~|~\Ga_{n}\wedge\Fa_n\right)=\Phi_{n+1}(m(\Xa_{n})
)^{\otimes (I_n)}\otimes \delta_{\textsl{z}_{n+1}}\otimes\Phi_{n+1}(m(\Xa_{n})
)^{\otimes (N-I_n)}.
$$
In addition, we have
\begin{equation}\label{f4}
\EE\left(m(\Xa_{n+1})(f)~|~\Ga_{n}\wedge\Fa_n\right)=\frac{m(\Xa_{n})(\Qa_{n+1}(f))}{m(\Xa_{n})(\Qa_{n+1}(1))}=\EE\left(m(\Xa_{n+1})(f)|\Ga_{n}\right).
\end{equation}
Therefore, using Lemma~\ref{L2}, for any non-negative functions $f,g\in \Ba(S_{n+1})$, s.t. $0<a\leq g\leq b<\infty$  for some finite constant $(a,b)$,
we have
$$
1-\frac{c_1}{N}~\leq \EE\left[
\frac{m(\Xa_{n+1})(f)}{m(\Xa_{n+1})(g)}~\Big|~\Ga_n
\right]~\frac{\EE\left[m(\Xa_{n+1})(g)~|~\Ga_n\right]}{\EE\left[m(\Xa_{n+1})(f)~|~\Ga_n\right]}\leq 
1+\frac{c_2}{N}
$$
as soon as $N>2(c_1+1)$, with
$$
c_1=\frac{b}{a}-1\quad\mbox{and}\quad c_2:= 2 c_1~(2c_1(c_1+1)+1)
$$

On the other hand, by (\ref{f4}) we have
$$
\frac{\EE\left[m(\Xa_{n+1})(f)~|~\Ga_n\right]}{\EE\left[m(\Xa_{n+1})(g)~|~\Ga_n\right]}=\frac{m(\Xa_{n})(\Qa_{n+1}(f))}{m(\Xa_{n})(\Qa_{n+1}(g))}.
$$
This yields
\begin{eqnarray*}
\left(1-\frac{c_1}{N}\right)\frac{m(\Xa_{n})(\Qa_{n+1}(f))}{m(\Xa_{n})(\Qa_{n+1}(g))}& \leq & \EE\left[
\frac{m(\Xa_{n+1})(f)}{m(\Xa_{n+1})(g)}~\Big|~\Ga_n
\right]\\ &\leq &
\left(1+\frac{c_2}{N}\right)\frac{m(\Xa_{n})(\Qa_{n+1}(f))}{m(\Xa_{n})(\Qa_{n+1}(g))}.
\end{eqnarray*}
Using (\ref{i2}) we conclude that
\begin{eqnarray*}
\left(1-\frac{c_1(\rho)}{N}\right)\frac{m(\Xa_{p})(\Qa_{p,n}(f))}{m(\Xa_{p})(\Qa_{p,n}(1))}& \leq &\EE\left[
\frac{m(\Xa_{p+1})(\Qa_{p+1,n}(f))}{m(\Xa_{p+1})(\Qa_{p+1,n}(1))}~\Big|~\Ga_p
\right]\\ 
& \leq &
\left(1+\frac{c_2(\rho)}{N}\right)\frac{m(\Xa_{p})(\Qa_{p,n}(f))}{m(\Xa_{p})(\Qa_{p,n}(1))}
\end{eqnarray*}
for any $p<n$, with
$$
c_1(\rho)=\frac{5}{2}~\rho^2-1\quad\mbox{and}\quad c_2(\rho):= 2 c_1(\rho)~(2c_1(\rho)(c_1(\rho)+1)+1)
$$

Iterating these estimates, we prove that
$$
\left(1-\frac{c_1(\rho)}{N}\right)^n\frac{\mu_0\Qa_{0,n}(f)}{\mu_0\Qa_{0,n}(1)}~\leq \EE\left[
m(\Xa_n)(f)
\right]~\leq 
\left(1+\frac{c_2(\rho)}{N}\right)^n\frac{\mu_0\Qa_{0,n}(f)}{\mu_0\Qa_{0,n}(1)}.
$$
The end of the proof of (\ref{f5}) is now easily completed.

Now we come to the proof of (\ref{f12}). 

{\bf Proof of (\ref{f12}) :}\label{pf12}

The proof of (\ref{f12}) is based on the transfer formula
\begin{equation}\label{f7}
m(x_{n-1})Q_{n}^{\prime}\HH_{n,m(x)}=m(x_{n-1})\HH_{n-1,m(x)}Q_n
\end{equation}
which is valid for any sequence $x=(x_k)_{k\geq 0}\in \prod_{k\geq 0}S^{\prime N}_k$.

We check this claim using the fact that
$$
 (m(x_{n-1})Q_{n}^{\prime})(d\textsl{x}_n)=m(x_{n-1})\left(H^{\prime}_{n}(\point,\textsl{x}_n)\right)~\nu^{\prime}_n(d\textsl{x}_n).
$$
This readily implies that
$$
\begin{array}{l}
( m(x_{n-1})Q_n^{\prime})\HH_{n,m(x)}(f)\\
\\
=\int m(x_{n-1})(d\textsl{y}_{n-1}) Q^{\prime}_{n}(\textsl{y}_{n-1},d\textsl{y}_{n})
\left\{\prod_{0\leq p<(n-1)}\frac{m(x_k)(d\textsl{y}_{k}) H^{\prime}_{k+1}(\textsl{y}_{k},\textsl{y}_{k+1})}{m(x_k)\left(H^{\prime}_{k+1}(\point,\textsl{y}_{k+1})\right)}\right\} f(\textsl{y}_0,\ldots,\textsl{y}_{n})\\
\\
=\int m(x_{n-1})(d\textsl{y}_{n-1}) 
\left\{\prod_{0\leq p<(n-1)}\frac{m(x_k)(d\textsl{y}_{k}) H^{\prime}_{k+1}(\textsl{y}_{k},\textsl{y}_{k+1})}{m(x_k)\left(H^{\prime}_{k+1}(\point,\textsl{y}_{k+1})\right)}\right\} Q_n(f)(\textsl{y}_0,\ldots,\textsl{y}_{n-1})\\
\\
=m(x_{n-1})\HH_{n-1,m(x)}(Q_n(f)).
\end{array}
$$

We let $\Qa^{\prime}_{n}$ be the integral operators defined as $\Qa_n$ by replacing $(S_n,Q_n)$ by $(S^{\prime}_n,Q^{\prime}_n)$; that is, we have that
\begin{equation}\label{f8}
\Qa^{\prime}_n:=Q^{\prime}_n\Ia^{\prime}_n \geq \left(1-\frac{1}{N}\right) Q^{\prime}_n 
\end{equation}
with
$$ \Ia_n^{\prime}(\textsl{x}^{\prime}_n,d\textsl{y}^{\prime}_n):=\frac{1}{N}~\delta_{\textsl{z}^{\prime}_n}(d\textsl{y}_n^{\prime})+\left(1-\frac{1}{N}\right)~\delta_{\textsl{x}^{\prime}_n}(d\textsl{y}^{\prime}_n).
$$
We denote by $\Qa^{\prime}_{p,n}$  the semigroup defined by the backward recursion 
$\Qa_{p-1,n}^{\prime}=\Qa_{p}^{\prime}\Qa_{p,n}^{\prime}$, 
Notice for any path sequence $\textsl{x}_p=(\textsl{x}^{\prime}_0,\ldots,\textsl{x}^{\prime}_p)\in S_p$ we have
\begin{equation}\label{f9}
\Qa_{p,n}(1)(\textsl{x}_p)=\Qa^{\prime}_{p,n}(1)(\textsl{x}_p^{\prime})=\EE\left(\prod_{p\leq k<n}G^{\prime}_k(X^{\prime}_{\textsl{z}^{\prime}\!,k})~\Big|~X^{\prime}_{\textsl{z}^{\prime}\!,p}=\textsl{x}_p^{\prime}\right)
\end{equation}
where $X^{\prime}_{\textsl{z}^{\prime},k}$ is a Markov chain with transitions $M^{\prime}_k\Ia_k^{\prime}$. This shows that the estimates (\ref{i2})
stated in Proposition~\ref{P1} are also satisfied if we replace $\Qa_{p,n}$ by the semigroup $\Qa^{\prime}_{p,n}$.

Now, let $\Ga_n^{\prime}=\sigma(\Xa^{\prime}_p,~p\leq n)$ be the natural $\sigma$-algebra filtration associated with the dual particle model $\Xa^{\prime}_{n}$.
By construction, for any $f\in\Ba(S^{\prime}_{n+1})$ we have
$$
\EE\left(m(\Xa_{n+1}^{\prime})(f)~|~\Ga_n^{\prime}\right)=m(\Xa_{n}^{\prime})(\Qa^{\prime}_{n+1}(f))/m(\Xa_{n}^{\prime})(\Qa^{\prime}_{n+1}(1)).
$$
Combining Lemma~\ref{L2} with (\ref{f7}) and (\ref{f8}), for any non-negative functions $f,g\in \Ba(S_{n+1})$ s.t. $0<a\leq g\leq b<\infty$   for some finite constant $(a,b)$, we prove that
$$
\begin{array}{l}
\EE\left[\frac{m(\Xa_{n+1}^{\prime})(\HH_{n+1,m(\Xa^{\prime})}(f))}{m(\Xa_{n+1}^{\prime})(g)}~\left\vert~\Ga_{n}~\right.\right]\\
\\
\geq \left(1-\frac{(b/a)-1}{N}\right)\left(1-\frac{1}{N}\right)~\frac{m(\Xa_{n}^{\prime})(\HH_{n,m(\Xa^{\prime})}(Q_{n+1}f))}{m(\Xa_{n}^{\prime})(\Qa^{\prime}_{n+1}(g))}
\geq \left(1-\frac{c}{N}\right)~\frac{m(\Xa_{n}^{\prime})(\HH_{n,m(\Xa^{\prime})}(Q_{n+1}f))}{m(\Xa_{n}^{\prime})(\Qa^{\prime}_{n+1}(g))}
 \end{array}
$$
for any $N\geq c:=(b/a)$. Combining (\ref{f9}) with (\ref{i2}) we readily check that
\begin{equation}\label{f10}
\EE\left[\frac{m(\Xa_{p}^{\prime})(\HH_{p,m(\Xa^{\prime})}(Q_{p,n}(f))}{m(\Xa_{p}^{\prime})(\Qa_{p,n}^{\prime}(1))}~\left\vert~\Ga_{p-1}~\right.\right]\geq 
\left(1-\frac{c}{N}\right)\frac{m(\Xa_{p-1}^{\prime})(\HH_{p-1,m(\Xa^{\prime})}(Q_{p-1,n}(f))}{m(\Xa_{p-1}^{\prime})(\Qa_{p-1,n}^{\prime}(1))}
\end{equation}
for any $N\geq c:=(5/2)\rho^2$. On the other hand, we have
$$
\begin{array}{l}
\EE\left(m(\Xa_{n}^{\prime})(\HH_{n,m(\Xa^{\prime})}(f)\right)\\
\\
\geq \left(1-\frac{1}{N}\right)\EE\left(\frac{m(\Xa_{n-1}^{\prime})Q_n^{\prime}(\HH_{n,m(\Xa^{\prime})}(f))}{m(\Xa_{n-1}^{\prime})(\Qa_{n}^{\prime}(1))}\right)= \left(1-\frac{1}{N}\right)~\EE\left[
\frac{m(\Xa_{n-1}^{\prime})(\HH_{n-1,m(\Xa^{\prime})}(Q_{n}(f))}{m(\Xa_{n-1}^{\prime})(\Qa_{n}^{\prime}(1))}
\right].
\end{array}
$$
Iterating the estimates (\ref{f10}), we conclude that
\begin{eqnarray*}
\EE\left(m(\Xa_{n}^{\prime})(\HH_{n,m(\Xa^{\prime})}(f)\right) &\geq& \left(1-\frac{1}{N}\right)\left(1-\frac{c}{N}\right)^{n}~\frac{\eta_0Q_{0,n}(f)}{\mu_0\Qa_{0,n}^{\prime}(1)}
\\ & \geq & \left(1-\frac{1+cn}{N}\right)~\frac{1}{\mu_0\Qa_{0,n}^{\prime}(1)}~\eta_{n}(f).
\end{eqnarray*}
Recalling that $\mu_0\Qa_{0,n}^{\prime}(1)=\mu_0\Qa_{0,n}(1)$ (see for instance (\ref{f9})), using (\ref{i3}) we check that
$$
\EE\left(m(\Xa_{n}^{\prime})(\HH_{n,m(\Xa^{\prime})}(f)\right)\geq\frac{1-\frac{1+cn}{N}}{1-\frac{1+n}{N}}~\eta_n(f)\geq \left(1-\frac{2+(c+1)n}{N}\right)~\eta_n(f).
$$

\appendix

\section{Main Technical Results}

We consider a collection $X=(X^1,\ldots,X^N)$ of $N\geq 1$ conditionally independent random variables w.r.t.~some $\sigma$-field $\Fa$,
 taking values in some measurable state space $(E,\Ea)$.
We set $[N]:=\{1,\ldots,N\}$, $[N]_i:=[N]-\{i\}$, and $\Fa_i:=\Fa\wedge \sigma(X^i)$, for any $i\in [N]$.
We consider the empirical measures
$$
m(X)=\frac{1}{N}\sum_{j\in [N]}\delta_{X^j}\quad\mbox{\rm and}\quad m_i(X):=\frac{1}{N-1}\sum_{j\in [N]_i}\delta_{X^j}~,\quad\mbox{\rm for any $i\in[N]$.}
$$

\begin{lem}\label{L1}
Let $f\in\Ba(E)$ be s.t.~$0<a\leq f\leq b<\infty$  for some finite constants $(a,b)$. Then, for any $i\in [N]$,
we have
\begin{equation}\label{f1}
 \left\vert\frac{\EE\left(m(X)(f)|\Fa_i\right)}{\EE\left(m(X)(f)|\Fa\right)}-1\right\vert\leq \frac{1}{N}~\left(\frac{b}{a}-1\right)
\end{equation}
and
\begin{equation}\label{f2}
1\leq \EE\left(m(X)(f)|\Fa\right)~\EE\left(\frac{1}{m(X)(f)}~\Big|~\Fa\right)\leq 1+\frac{1}{N}~\frac{b}{a}~\left(\frac{b}{a}-1\right)^2.
\end{equation}
\end{lem}
\proof
Formula (\ref{f1}) comes from the fact that
\begin{eqnarray*}
\EE\left(m(X)(f)~|~\Fa_i\right)&=&\frac{1}{N}~f(X^i)+\left(1-\frac{1}{N}\right)~\EE\left(m_i(X)(f)|\Fa\right)\\
&=&\EE\left(m(X)(f)|\Fa\right)+\frac{1}{N}~\left[f(X^i)-\EE\left(f(X^i)|\Fa\right)\right].
\end{eqnarray*}
The l.h.s.~of formula (\ref{f2}) is a direct consequence of Jensen's inequality. To check the r.h.s. estimate, we use the decomposition
\begin{eqnarray*}
\frac{1}{m(X)(\overline{f})}&=&1+\left(1-m(X)(\overline{f})\right)+\frac{\left(m(X)(\overline{f})-1\right)^2}{m(X)(\overline{f})}\\
&\leq& 1+\left(1-m(X)(\overline{f})\right)+\frac{b}{a}~\left(m(X)(\overline{f})-1\right)^2\quad\mbox{\rm with}\quad
\overline{f}:=\frac{f}{\EE\left(m(X)(f)|\Fa\right)}.
\end{eqnarray*}

Taking the expectation, we find that
\begin{eqnarray*}
\EE\left(\frac{1}{m(X)(\overline{f})}~|~\Fa\right)&\leq& 1+\frac{1}{N}\frac{b}{a}~\frac{1}{N}\sum_{i\in[N]}\EE\left[(\overline{f}(X^i)-\EE(\overline{f}(X^i)|\Fa))^2|\Fa\right]\\
&\leq& 1+\frac{1}{N}~\frac{b}{a}~\frac{1}{a^2}~\left(b-a\right)^2=1+\frac{1}{N}~\frac{b}{a}~\left(\frac{b}{a}-1\right)^2.
\end{eqnarray*}
This ends the proof of the lemma.
\cqfd

\begin{prop}\label{p1}
For any $f,g\in\Ba(E)$ be s.t.~$0<a\leq g\leq b<\infty$, $f:E\rightarrow\mathbb{R}^+$  for some finite constants $(a,b)$, we have the estimates
\begin{equation}\label{f18}
\left(1-\frac{c_1}{N}\right)
~\frac{\EE\left(m(X)(f)|\Fa\right)}{\EE\left(m(X)(g)|\Fa\right)}
~\leq \EE\left(
\frac{m(X)(f)}{m(X)(g)}
~\Big|~\Fa\right)\leq 
\frac{\EE\left(m(X)(f)|\Fa\right)}{\EE\left(m(X)(g)|\Fa\right)}~\left(1+\frac{c_2}{N}\right)
\end{equation}
as soon as $N>2(c_1+1)$, with
$$
c_1=\frac{b}{a}-1\quad\mbox{and}\quad c_2:= 2 c_1~(2c_1(c_1+1)+1).
$$
\end{prop}
\proof
Using (\ref{f1}),we have
$$
\left(1-\frac{c}{N}\right)~ \EE\left(m(X)(g)|\Fa\right)\leq \EE\left(m(X)(g)~|~\Fa_i\right)\leq \left(1+\frac{c}{N}\right)~ \EE\left(m(X)(g)|\Fa\right).
$$
with $c=\left(\frac{b}{a}-1\right)$.
We also observe that
$
m(X)(g)=m_i(X)\left(g_i\right)$,  with 
$$
 a\leq g_i:=\frac{1}{N}~g(X^i)+\left(1-\frac{1}{N}\right)~g\leq b.
$$
Thus, applying (\ref{f2})  to $X=(X^j)_{j\in [N]_i}$, $\Fa=\Fa_i$, and $f=g_i$ we prove that
$$
\frac{1}{1+\frac{c}{N}}~\frac{1}{ \EE\left(m(X)(g)|\Fa\right)}
\leq\frac{1}{ \EE\left(m_i(X)(g_i)|\Fa_i\right)}\leq  \EE\left(1/m_i(X)(g_i)|\Fa_i\right)=\EE\left(1/m(X)(g)|\Fa_i\right)$$
and for any $N>(c\vee 1)$
\begin{eqnarray*}
\EE\left(1/m(X)(g)|\Fa_i\right)= \EE\left(1/m_i(X)(g_i)|\Fa_i\right)&\leq& \frac{1}{ \EE\left(m_i(X)(g_i)|\Fa_i\right)}
\left(1+\frac{1}{N-1}~(c+1)~c^2\right)\\
&\leq&~\frac{1}{ \EE\left(m(X)(g)|\Fa\right)}\times \frac{\left(1+\frac{2}{N}~(c+1)~c^2\right)}{\left(1-\frac{c}{N}\right)}.
\end{eqnarray*}

For any positive constants $c_1$ and $2c_2<N$, we have
\begin{equation}\label{f25}
\frac{1+c_1/N}{1-c_2/N}\leq 1+2(c_1+c_2)/N\quad\mbox{\rm and}\quad \frac{1}{1+c_1/N}\geq 1-c_1/N.
\end{equation}
We check these claims using the fact that $\frac{1+x}{1-y}-1=\frac{x+y}{1-y}\leq 2(x+y)$, for any $x+y>0$ and $2y\leq 1$.
This yields the estimates
$$
\left(1-\frac{c}{N}\right)~\frac{1}{ \EE\left(m(X)(g)|\Fa\right)}\leq  \EE\left(\frac{1}{m(X)(g)}\Big|\Fa_i\right)$$
and
$$
  \EE\left(\frac{1}{m(X)(g)}\Big|\Fa_i\right)\leq \frac{1}{ \EE\left(m(X)(g)|\Fa\right)}\times\left(
1+\frac{2}{N}~\left[  2(c+1)~c^2+c\right]\right)
$$
as soon as $N>2(c+1)$.

On the other hand, we have
$$
\EE\left(
\frac{m(X)(f)}{m(X)(g)}
\Big|\Fa\right)=\frac{1}{N}\sum_{i\in[N]}~\EE\left(\EE\left(\frac{1}{m(X)(g)}~\Big|~\Fa_i\right)~f(X^i)~\Big|~\Fa\right).
$$
The end of the proof of (\ref{f18}) is now easily completed.
\cqfd

\section{Tensor Product Measures}\label{tensor-section}

We let $X=\left(X^i\right)_{1\leq i\leq N}$ be a sequence of independent random variables on some state space $E$. For any $q<N$ we set
$$
m(X)^{\odot q}=\frac{1}{(N)_q}\sum_{c\in I_q^N}\delta_{\left(X^{c(1)},\ldots,X^{c(q)}\right)}
$$
We recall that $I_q^N$ stands the set of all $(N)_q=\frac{(N)!}{(N-q)!}$ multi-indexes 
 $c=(c(1),\ldots,c(q))\in\{1,\ldots,N\}^q$ with different values, or equivalently one to one mappings from $[q]:=\{1,\ldots,q\}$ into
 $[N]:=\{1,\ldots,N\}$. Before we discuss first order estimates of tensor product measures, we provide some
useful  combinatorial estimates. Firstly, we observe that
\begin{equation}\label{f26}
1-\frac{(q-1)^2}{N}\leq \left(1-\frac{q-1}{N}\right)^{q-1}\leq \prod_{1\leq k<q}\left(1-\frac{k}{N}\right)=\frac{(N)_q}{N^q}\leq 1.
 \end{equation}
The l.h.s.~estimate comes from the fact that $1-x^{p}=(1-x)~\sum_{0\leq l<p}x^p\leq p(1-x)$, for any $0\leq x<1$ (choosing 
$p=(q-1)$ and $x=1-\frac{q-1}{N}$, we find that $1-\left(1-\frac{q-1}{N}\right)^{q-1}\leq \frac{(q-1)^2}{N}$).
$$
\frac{1}{N^q}\left[(N)_q-(N-q)_{q}\right] 
$$
$$
=\left[\prod_{1\leq i\leq q}\left(1-\frac{i-1}{N}\right)-\prod_{1\leq i\leq q}\left(1-\frac{q+(i-1)}{N}\right)\right]
$$
$$
=\sum_{1\leq k\leq q}\sum_{1\leq i_1<\ldots<i_k\leq q}\left\{\prod_{0\leq k<q}\left[\left(1-\frac{i_l-1}{N}\right)-\left(1-\frac{q+(i_l-1)}{N}\right)\right]\right\}
$$
$$
\hskip6cm\times\prod_{j\in\{1,\ldots,q\}-\{i_1,\ldots,i_k\}}\left(1-\frac{q+(j-1)}{N}\right)
$$
$$
\leq \sum_{1\leq k\leq q}\left(\begin{array}{c}
q\\
k
\end{array}
\right)~\left(\frac{q}{N}\right)^k~\left(1-\frac{q}{N}\right)^{q-k}=1-\left(1-\frac{q}{N}\right)^q\leq \frac{q^2}{N}.
$$
We conclude that
$$
\frac{1}{(N)_q}~\left[(N)_q-(N-q)_{q}\right]\leq \frac{N^q}{(N)_q}~\frac{q^2}{N}\leq \frac{1}{1-\frac{(q-1)^2}{N}}~~\frac{q^2}{N}\leq \frac{2q^2}{N}
$$
and
\begin{equation}\label{f21}
(N-q)_{q}=\frac{(N)_{2q}}{(N)_q}\Longrightarrow~
\frac{1}{(N)_q}~\left[(N)_q-(N-q)_{q}\right]=\frac{1}{(N)_q^2}~\left[(N)_q^2-(N)_{2q}\right]\leq \frac{2q^2}{N}
\end{equation}
as soon as $N\geq 2(q-1)^2$.

 Given a subset $J\subset [N]$ with cardinality $\#(J)=q$, we consider the empirical measure
 associated with multi-indices with different values and avoiding the set $J$; that is, we have
 $$
 m_J(X)^{\odot q}=\frac{1}{(N-q)_q}\sum_{c\in I_{q,J}^N}\delta_{\left(X^{c(1)},\ldots,X^{c(q)}\right)}
 $$
 where $I_{q,J}^N$ stands for the set of all one to one mappings from $[q]:=\{1,\ldots,q\}$ into
 $[N]-J$. We also set $\Fa_J:=\sigma(X^j,~j\in J)$ the $\sigma$-field generated by the random variables $X^j$ indexed by
 $j\in J$.

 The link between these measures and tensor product measures 
 is expressed in terms of the  Markov transitions
 $C_{c}$ indexed by the set of  mappings $c$ from $[q]$ into itself.  $C_{c}$  is defined for any $x=(\textsl{x}^1,\ldots,\textsl{x}^q)\in E^q$
and any function $f$ on $\Ba(E^q)$,
 by 
 $$
 C_{c}(f)(x)=f(x^c)\quad\mbox{\rm with}\quad
 x^c:=\left(\textsl{x}^{c(1)},\ldots,\textsl{x}^{c(q)}\right).
$$

 We emphasize that the tensor product measures discussed above are symmetry-invariant by construction. In the further development of this section,
 it is assumed without restriction that these measures act on symmetric functions $F$; that is
 $
 F=\frac{1}{q!}\sum_{\sigma\in\Ga_q} C_{\sigma}(F)
 $,
where $\Ga_q$ stands for the symmetric group of all permutations of $[q]$.
 The connection between these measures 
 is described in the following technical lemma taken from~\cite{dpr-2009}.

\begin{lem}\label{L3}
For any $q<N$ we have the formula
 $$
m(X)^{\otimes q}= m(X)^{\odot q}\Ca
\quad\mbox{with}\quad
 \Ca=\frac{1}{(N)^q}\sum_{c\in [q]^{[q]}}~\frac{(N)_{\vert c\vert}}{(q)_{\vert c\vert}}~C_c
 $$
where $\vert c\vert$ for the cardinality of the set $c([q])$, and $(m)_p=m!/(m-p)!$ stands for the number of one to one mappings
from $[p]$ into $[m]$.
\end{lem}

The next Lemma is an extension of the binomial-type falling factorial formula (a.k.a.~the Vandermonde convolution) 
to tensor product measures.

\begin{lem}\label{L4}
For any subset $J\subset [N]$ with cardinality $\#(J)=q<N$ and any
function $f$ on $E^q$, we have
\begin{eqnarray}\label{f19}
m(X)^{\odot q}(f)&=&m_J(X)^{\odot q}(f_{J})\nonumber\\
&=&\frac{1}{(N)_q}
\sum_{0\leq k\leq q}
\left(\begin{array}{c}
q\\
k
\end{array}
\right)
(q)_k~(N-q)_{q-k}~ m_J(X)^{\odot q}(f_{k,J})
\end{eqnarray}
with the random functions
\begin{equation}\label{f23}
f_{J}=\frac{1}{(N)_q}
\sum_{0\leq k\leq q}
\left(\begin{array}{c}
q\\
k
\end{array}
\right)
(q)_k~(N-q)_{q-k}~f_{k,J}
\end{equation}
and
$$
f_{k,J}\left(x^{1},\ldots,x^q\right):=\frac{1}{(q)_k}\sum_{a\in I_{k,J}^{(1,N)}}f\left(X^a,(x^{k+1},\ldots,x^q)\right).
$$
In the above display, $I_{k}^{J}$ stands for the set of all $(q)_k$ one to one mappings from $\{1,\ldots,k\}$ into
 $J$.
\end{lem}

\proof
We have
$$
\sum_{c\in I_q^N}f\left(X^c\right)=\sum_{0\leq k\leq q}
\left(\begin{array}{c}
q\\
k
\end{array}
\right)
\sum_{a\in I_{k,J}^{(1,N)}}\sum_{b\in I_{k,J}^{(2,N)}}f\left(X^a,X^b\right)
$$
where $I_{k,J}^{(1,N)}$ stands for the set of all $(q)_k$ one to one mappings from $\{1,\ldots,k\}$ into
 $J$; and $I_{k,J}^{(2,N)}$  stands for the set of all $(N-q)_{q-k}$ one to one mappings from $\{k+1,\ldots,k+(q-k)\}$ into
 $[N]-J$. 
 
 On the other hand, we have
$$
\sum_{a\in I_{k,J}^{(1,N)}}\sum_{b\in I_{k,J}^{(2,N)}}f\left(X^a,X^b\right)=(q)_k~\sum_{b\in I_{k,J}^{(2,N)}}\overline{f}_{k,J}\left(X^b\right)
$$ 
 with the random functions
 $$
 \overline{f}_{k,J}\left(x^{k+1},\ldots,x^q\right):=\frac{1}{(q)_k}\sum_{a\in I_{k,J}^{(1,N)}}f\left(X^a,(x^{k+1},\ldots,x^q)\right).
 $$
 Notice that for any function $g$ on $E^{q-k}$, we have
$$
\frac{1}{(N-q)_q}\sum_{c\in I_{q,J}^N}g\left(X^{c(k+1)},\ldots,X^{c(q)}\right)=\frac{1}{(N-q)_{q-k}}\sum_{c\in I_{q,J}^{(2,N)}}g\left(X^{c(k+1)},\ldots,X^{c(q)}\right).
$$
This yields
$$
 m_J(X)^{\odot q}(f_{k,J})=\frac{1}{(N-q)_{q-k}}\sum_{c\in I_{q,J}^{(2,N)}} \overline{f}_{k,J}\left(X^{c(k+1)},\ldots,X^{c(q)}\right)
$$
and
$$
\sum_{a\in I_{k,J}^{(1,N)}}\sum_{b\in I_{k,J}^{(2,N)}}f\left(X^a,X^b\right)=(q)_k~(N-q)_{q-k}~ m_J(X)^{\odot q}(f_{k,J}).
$$
This clearly ends the proof of (\ref{f19}).
\cqfd

\begin{lem}\label{L5}
Let $f\in\Ba(E^q)$ be s.t.~$0<a\leq f\leq b<\infty$  for some finite constants $(a,b)$. 
Then, for any given a subset $J\subset [N]$ with cardinality $\#(J)=q$, and any $N\geq 2q^2$
we have
\begin{equation}\label{f20}
\left\vert\frac{\EE\left(m(X)^{\odot q}(f)|\Fa_{J}\right)}{\EE\left(m(X)^{\odot q}(f)\right)}-1\right\vert\leq \frac{2q^2}{N}~\left(\frac{b}{a}-1\right)
\end{equation}
and
\begin{equation}\label{f22}
1\leq \EE\left(m(X)^{\odot q}(f)\right)
\EE\left(\frac{1}{m(X)^{\odot q}(f)}\right)
\leq1+\frac{2q^2}{N}~\frac{b}{a}~\left(\frac{b}{a}-1\right)^2.
\end{equation}
\end{lem}
\proof
Combining (\ref{f19}), 
with the fact that 
 $$f_{0,J}=f\quad\mbox{\rm and}\quad \EE\left(m_J(X)^{\odot q}(f)|\Fa_{J}\right)=\EE\left(m(X)^{\odot q}(f)\right)
 $$
we prove that
$$
\begin{array}{l}
\EE\left(m(X)^{\odot q}(f)|\Fa_{J}\right)\\
\\
=\EE\left(m_J(X)^{\odot q}(f_J)|\Fa_{J}\right)\\
\\
=\EE\left(m(X)^{\odot q}(f)\right)
+\frac{1}{(N)_q}
\sum_{1\leq k\leq q}
\left(\begin{array}{c}
q\\
k
\end{array}
\right)
(q)_k~(N-q)_{q-k}~\\
\\
\hskip5cm\times \left[\EE\left(m_J(X)^{\odot q}(f_{k,J})|\Fa_{J}\right)-\EE\left(m_J(X)^{\odot q}(f)|\Fa_{J}\right)\right].
\end{array}
$$
This implies that
\begin{eqnarray*}
 \left\vert\frac{\EE\left(m(X)^{\odot q}(f)|\Fa_{J}\right)}{\EE\left(m(X)^{\odot q}(f)\right)}-1\right\vert&\leq& 
 \frac{1}{(N)_q}
\sum_{1\leq k\leq q}
\left(\begin{array}{c}
q\\
k
\end{array}
\right)
(q)_k~(N-q)_{q-k} \left(\frac{b-a}{a}\right)\\
&=&\left(1-\frac{(N-q)_q}{(N)_q}\right) \left(\frac{b}{a}-1\right)\leq \frac{2q^2}{N}~\left(\frac{b}{a}-1\right).
\end{eqnarray*}
In the last assertion, we have used (\ref{f21}) and Vandermonde convolution formula
$$
 \frac{1}{(N)_q}
\sum_{0\leq k\leq q}
\left(\begin{array}{c}
q\\
k
\end{array}
\right)
(q)_k~(N-q)_{q-k} =1.
$$
This ends the proof of (\ref{f20}). The proof of (\ref{f22}) follows the same line of arguments as the proof of (\ref{f2}). Indeed, if we set
$\overline{f}:=\frac{f}{\EE\left(m(X)^{\odot q}(f)\right)}$ the r.h.s. estimates comes from the fact that
$$
\EE\left(\left[m(X)^{\odot q}(\overline{f}-1)\right]^2\right)=\frac{1}{(N)_q^2}\sum_{(a,b)\in (I^N_q)^2}\EE\left((\overline{f}(X^a)-1)(\overline{f}(X^b)-1)\right).
$$
On the other hand, we have
$$
a([q])\cap b([q])=\emptyset\Longrightarrow
\EE\left((\overline{f}(X^a)-1)(\overline{f}(X^b)-1)\right)=0
$$
and we have $(N)_{2q}$ possible pairs $(a,b)$ satisfying this property.
Using (\ref{f21}), this implies that
$$
\EE\left(\left[m(X)^{\odot q}(\overline{f}-1)\right]^2\right)\leq \frac{1}{(N)_q^2}\left[(N)_q^2-(N)_{2q}\right]~(b-a)^2\leq \frac{2q^2}{N}~(b-a)^2
$$
and therefore
$$
\EE\left(\frac{1}{m(X)^{\odot q}(\overline{f})}\right)\leq 1+\frac{b}{a}~\EE\left(\left[m(X)^{\odot q}(\overline{f}-1)\right]^2\right) \leq1+\frac{2q^2}{N}~\frac{b}{a}~\left(\frac{b}{a}-1\right)^2.
$$
This ends the proof of the lemma.
\cqfd

\begin{prop}\label{P3}
For any $f,g\in\Ba(E^q)$ s.t.~$0<a\leq g\leq b<\infty$ for some finite constants $(a,b)$, and $f:E^q\rightarrow\mathbb{R}^+$
s.t. $\EE\left(f(X^c)\right)<\infty$ for any $c\in I^N_q$, we have the estimates
$$
\left(1-\frac{c_1(q)}{N}\right)~\frac{\EE\left(m(X)^{\odot q}(f)\right)}{\EE\left(m(X)^{\odot q}(g)\right)}\leq 
\EE\left(
\frac{m(X)^{\odot q}(f)}{m(X)^{\odot q}(g)}
\right)
\leq 
 \frac{\EE\left(m(X)^{\odot q}(f)\right)}{\EE\left(m(X)^{\odot q}(g)\right)}~\left(1+\frac{c_2(q)}{N}\right)
$$
with
$$
c_1(q):=2q^2c\quad\mbox{and}\quad
c_2(q):=2c_1(q)~\left(1+2c(c+1)\right)\quad\mbox{with}\quad c=\left(\frac{b}{a}-1\right).
$$
\end{prop}

\proof
We use the decomposition
$$
\EE\left(
\frac{m(X)^{\odot q}(f)}{m(X)^{\odot q}(g)}
\right)=\frac{1}{(N)_q}\sum_{\alpha\in I^N_q}\EE\left(\EE\left(\frac{1}{m(X)^{\odot q}(g)}~|~\Fa_{J_{\alpha}}\right)~f(X^{\alpha})\right)
$$
with $J_{\alpha}=\alpha([q])=\{\alpha(1),\ldots,\alpha(q)\}$. By construction, we have
$$
m(X)^{\odot q}(g)=m_{J_{\alpha}}(X)^{\odot q}(g_{J_{\alpha}})
$$
with the $[a,b]$-valued random function $g_{J_{\alpha}}$ defined as in (\ref{f23}) by replacing $f$ by $g$ and $J$ by $J_a$. By Lemma~\ref{L5}, we have
\begin{eqnarray*}
\EE\left(\frac{1}{m_{J_{\alpha}}(X)^{\odot q}(g_{J_{\alpha}})}~|~\Fa_{J_{\alpha}}\right)&\geq &\frac{1}{\EE\left(m_{J_{\alpha}}(X)^{\odot q}(g_{J_{\alpha}})~|~\Fa_{J_{\alpha}}\right)}=
\frac{1}{\EE\left(m(X)^{\odot q}(g)~|~\Fa_{J_{\alpha}}\right)}\\
&\geq &
\frac{1}{1+\frac{2cq^2}{N}}~\frac{1}{\EE\left(m(X)^{\odot q}(g)\right)}\\
&\geq &
\left(1-\frac{2cq^2}{N}\right)~\frac{1}{\EE\left(m(X)^{\odot q}(g)\right)}
\end{eqnarray*}
and
\begin{eqnarray*}
 \EE\left(\frac{1}{m_{J_{\alpha}}(X)^{\odot q}(g_{J_{\alpha}})}~|~\Fa_{J_{\alpha}}\right)&\leq& 
\frac{1}{\EE\left(m_{J_{\alpha}}(X)^{\odot q}(g_{J_{\alpha}})~|~\Fa_{J_{\alpha}}\right)}~\left(1+\frac{2q^2}{N-q}~(c+1)~c^2\right)\\
&\leq & \frac{1}{\EE\left(m(X)^{\odot q}(g)\right)}~\frac{1+\frac{4q^2}{N}~(c+1)~c^2}{1-\frac{2cq^2}{N}}
\end{eqnarray*}
as soon as $N\geq 2q^2c$ with $c=\left(\frac{b}{a}-1\right)$. In the last estimate we have used the fact that
$\frac{1}{N-q}=\frac{1}{N}~\frac{1}{1-q/N}\leq \frac{2}{N}$ for any $N\geq 2q$. Using (\ref{f25}), we conclude that
$$
 \EE\left(\frac{1}{m_{J_{\alpha}}(X)^{\odot q}(g_{J_{\alpha}})}~|~\Fa_{J_{\alpha}}\right)\leq 
 \frac{1}{\EE\left(m(X)^{\odot q}(g)\right)}~\left(1+\frac{4cq^2}{N}\left(1+2c(c+1)\right)\right)
$$
for any $N>4cq^2$. This ends the proof of the proposition.

\cqfd

\section{Propagation of Chaos Estimates}\label{pce-section}
Our next objective is to prove the propagation of chaos estimates discussed in Section~\ref{propagation-chaos-sec}. 
By construction, we have
\begin{eqnarray*}
\EE\left(m(\xi_n)^{\odot q}(f)~|~\xi_{n-1}\right)&=&\frac{m(\xi_{n-1})^{\otimes q}(Q_n^{\otimes q}(f))}{m(\xi_{n-1})^{\otimes q}(Q_n^{\otimes q}(1))}
\\&=&\frac{m(\xi_{n-1})^{\odot q}(\Ca Q_n^{\otimes q}(f))}{m(\xi_{n-1})^{\odot q}(\Ca Q_n^{\otimes q}(1))}=\frac{m(\xi_{n-1})^{\odot q}(Q_n^{(q)}(f))}{m(\xi_{n-1})^{\odot q}(Q_n^{(q)}(1))}
\end{eqnarray*}
with $Q_n^{(q)}=\Ca Q_n^{\otimes q}$, and the coalescent Markov transition $\Ca$ introduced in Lemma~\ref{L3}.
We denote by $$
\forall 0\leq p\leq n\qquad Q_{p,n}^{(q)}=Q_{p+1}^{(q)}Q_{p+1,n}^{(q)}$$ 
the corresponding semigroup, with the convention $Q_{n,n}^{(q)}=Id$.

On the other hand, we have
\begin{eqnarray*}
\frac{\EE\left(m(\xi_{n-1})^{\odot q}(Q_n^{(q)}(f))~|~\xi_{n-2}\right)}{\EE\left(m(\xi_{n-1})^{\odot q}(Q_n^{(q)}(1))~|~\xi_{n-2}\right)}&=&
\frac{m(\xi_{n-2})^{\otimes q}(Q_{n-1}^{\otimes q}Q_n^{(q)}(f))}{m(\xi_{n-2})^{\otimes q}(Q_{n-1}^{\otimes q}Q_n^{(q)}(1))}\\
&=&
\frac{m(\xi_{n-2})^{\odot q}(Q_{n-2,n}^{(q)}(f))}{m(\xi_{n-2})^{\odot q}(Q_{n-2,n}^{(q)}(1))}.
\end{eqnarray*}

We will seek to sequentially
apply Proposition~\ref{P3}; first we need to estimate the functions $Q_{p,n}^{(q)}(1)$. We assume that (A\ref{hyp:a1}-\ref{hyp:a2}) 
are satisfied.

We rewrite $\Ca$ in terms of the coalescence degree
$$
 \Ca=\sum_{0\leq l<q}\frac{(N)_{q-l}}{N^{q-l}}~\frac{1}{N^{l}}~S(q,q-l)~\CC_l
$$
with the uniform coalescent transition
$$
\CC_l:=\frac{1}{S(q,q-l)(q)_{q-l}}~\sum_{c\in [q]^{[q]},~\vert c\vert=q-l} C_c.
$$
In the above display, $S(q,k)$ stands for the Stirling number of the second kind.% (the number of ways to partition a set of $q$ objects into $k$ non-empty subsets).
We recall from~\cite{rennie} that
$$
S(q,q-l)\leq \frac{1}{2}~\left(\begin{array}{c}q\\l\end{array}\right)~(q-l)^l\leq \left(\begin{array}{c}q\\l\end{array}\right)~q^l.
$$
By construction, for any sequence of coalescent indices $0\leq l_i<q$, with $p<i\leq n$ we have
$$
 (1/\rho)^{q+\sum_{p<i\leq n}l_i}\leq \left(\CC_{l_{p+1}}Q_{p+1}^{\otimes q}\right)\left(\CC_{l_{p+2}}Q_{p+2}^{\otimes q}\right)\ldots \left(\CC_{l_{n}}Q_{n}^{\otimes q}\right)(1)\leq \rho^{q+\sum_{p<i\leq n}l_i}
$$

Using (\ref{f26}), we also notice that
$$
1-\frac{q^2(n-p)}{N}\leq \left(1-\frac{q^2}{N}\right)^{n-p}
\leq
\prod_{p<i\leq n} \frac{(N)_{q-l_i}}{N^{q-l_i}}\leq 1.
$$
Using these estimates, we prove that
$$
\left(1-\frac{q^2(n-p)}{N}\right)~Q^{\otimes q}_{p,n}(f)
\leq
\left(\frac{(N)_q}{N^q}\right)^{n-p}
Q^{\otimes q}_{p,n}(f)\leq Q^{(q)}_{p,n}(f)
$$
and
\begin{eqnarray*}
Q^{(q)}_{p,n}(f)&\leq& \sum_{0\leq l_{p+1},\ldots,l_n<q}\left[\prod_{p<i\leq n} \left(\begin{array}{c}q\\l_i\end{array}\right)~\left(\frac{q}{N}\right)^{l_i}\right]
\left(\CC_{l_{p+1}}Q_{p+1}^{\otimes q}\right)\ldots \left(\CC_{l_{n}}Q_{n}^{\otimes q}\right)(f)\\
&\leq & Q^{\otimes q}_{p,n}(f)+\rho^q~\sum_{1\leq s<(n-p)q} \sum_{l_{p+1}+\ldots+l_n=s}\left[\prod_{p<i\leq n} \left(\begin{array}{c}q\\l_i\end{array}\right)~\left(\frac{q\rho }{N}\right)^{l_i}\right]\\
&=& Q^{\otimes q}_{p,n}(f)+\rho^q~\left[\left(1+\frac{q\rho}{N}\right)^{q(n-p)}-1\right]\leq Q^{\otimes q}_{p,n}(f)+\rho^{q+1}~eq^2~\frac{(n-p)}{N}
\end{eqnarray*}
for any non-negative function $f$ on $S_n^q$ s.t $\|f\|\leq 1$, and any $N\geq  e q^2(n-p)\rho$.
For non-negative tensor product functions
$f=g^{\otimes q}$, we also have
$$
\begin{array}{l}
Q_n^{\otimes q}(f)=Q_n(1)^q \left(\frac{Q_n(f)}{Q_n(1)}\right)^q\leq Q_n(1)^{q-1}~Q_n(g^q)=Q_n^{\otimes q}(1^{\otimes (q-1)}\otimes g^q)
\\
\\
\Longrightarrow \left(\CC_{l_{p+1}}Q_{p+1}^{\otimes q}\right)\ldots \left(\CC_{l_{n}}Q_{n}^{\otimes q}\right)(f) \leq \rho^{q+\sum_{p<i\leq n}l_i}\times  Q_{p,n}(g^q).
\end{array}
$$
Summarizing, we have proved the following estimates:
\begin{lem}
Assume (A\ref{hyp:a1}-\ref{hyp:a2}). Then, for any 
$f:S_n^q\rightarrow[0,1]$   and for  any $N\geq  e q^2(n-p)\rho$, we have
$$
\left(1-\frac{q^2(n-p)}{N}\right)~Q^{\otimes q}_{p,n}(f)\leq Q^{(q)}_{p,n}(f)\leq Q^{\otimes q}_{p,n}(f)+\rho^{q+1}~eq^2~\frac{(n-p)}{N}
$$
In addition, for tensor product functions
$f=g^{\otimes q}$ we have
$$
Q^{(q)}_{p,n}(g^{\otimes q})\leq Q^{\otimes q}_{p,n}(g^{\otimes q})+\rho^{q+1}~eq^2~\frac{(n-p)}{N}~Q_{p,n}(g^q).
$$
\end{lem}

This lemma readily implies that
$$
N\geq  e q^2n\rho\Rightarrow (2\rho^{q})^{-1}\leq \left(1-\frac{q^2n}{N}\right)~\rho^{-q}~\leq Q^{(q)}_{p,n}(1)\leq \rho^q\left(1+\rho~eq^2~\frac{n}{N}\right)\leq 2\rho^q
$$
as well as
$$
\left(1-\frac{q^2n}{N}\right)~\eta^{\otimes q}_{n}(f)\leq \eta_0^{\otimes q}Q^{(q)}_{0,n}(f)\leq \eta^{\otimes q}_{n}(f)+\rho^{q+1}~eq^2~\frac{n}{N}
$$
and
$$
\left(1-\frac{q^2n}{N}\right)~(\eta_{n}(g))^q\leq \eta_0^{\otimes q}Q^{(q)}_{0,n}(g^{\otimes q})\leq(\eta_{n}(g))^q+\rho^{q+1}~eq^2~\frac{n}{N}~\eta_n(g^q).
$$

If we set
$$
\gamma^{(q)}_n(f):= \eta_0^{\otimes q}Q^{(q)}_{0,n}(f)\quad\mbox{\rm and}\quad \eta^{(q)}_n(f):=\gamma^{(q)}_n(f)/\gamma^{(q)}_n(1)
$$
then, sequentially iterating Proposition~\ref{P3}, we prove the following estimates:
\begin{prop}
Assume (A\ref{hyp:a1}-\ref{hyp:a2}). Then, for any function $f:S_n^q\rightarrow\mathbb{R}^+$ with $\eta_n^{(q)}(f)<\infty$ we have
$$
\left(1-a(\rho)\frac{n}{N}\right)
 \eta^{(q)}_n(f)
\leq \EE\left(m(\xi_n)^{\odot q}(f)\right)\leq  \eta^{(q)}_n(f)\left(1+a(\rho)\frac{n}{N}\right)
$$
for some finite constants $a(\rho)$ and as soon as $N\geq n a(\rho)$. 
\end{prop}

The end of the proof of the propagation of chaos estimates discussed
in Section~\ref{propagation-chaos-sec} for $q$-tensor product empirical measures now follows the same line of arguments as the ones we used in the proof of Theorem~\ref{T2} for particle density profiles (i.e. for $q=1$).


\begin{thebibliography}{99}    


\bibitem{adh-2010}
C. Andrieu, A. Doucet, and  R. Holenstein. Particle Markov chain Monte Carlo methods
(with discussion). \emph{J. R. Statist. Soc. B}, 72, 1-269 (2010).

\bibitem{alv-2013}
C. Andrieu, A. Lee, and M. Vihola. Uniform ergodicity of the iterated conditional SMC and geometric
ergodicity of particle Gibbs samplers. arXiv:1312.6432v1 (2013).


\bibitem {caffarel}R. Assaraf and M. Caffarel. A pedagogical introduction to
quantum Monte Carlo. In \emph{Mathematical Models and Methods for Ab Initio
Quantum Chemistry}, Lecture Notes in Chemistry, eds. M. Defranceschi and C. Le
Bris, Springer (2000).

\bibitem{cances-tony}
E. Canc\`es,  B. Jourdain, and T. Leli\`evre. Quantum Monte Carlo simulations of fermions. A mathematical analysis of the fixed-node approximation, \emph{Math. Models Methods Appl. Sci.}, 16, 1403-1440, (2006).

\bibitem{cappe-moulines}
O. Capp\'e, E. Moulines, and T. Ryd\`en. Inference in Hidden Markov Models,
Springer-Verlag (2005).

\bibitem{cdho}
R. Carmona, P. Del Moral, P. Hu, and N. Oudjane An introduction to particle methods in finance 
in Numerical Methods in Finance Springer New York, Series: Proceeding in Mathematics
Vol. 12 (2012). 

\bibitem{cdmg-11}
F. Cerou, P. Del Moral, and A. Guyader.
A non-asymptotic variance theorem for unnormalized Feynman-Kac particle models
\emph{Ann. Inst. H. Poincar\'e Probab. Statist.}, 47,  629-649 (2011). 

\bibitem{cl-2013}
H. P. Chan and T. L. Lai.
A general theory of particle filters in hidden Markov models and some applications. \emph{Ann. Statist.}, 
41, 2877-2904 (2013).

\bibitem{d-2004}
P. {Del Moral}.
\newblock \href{http://www.math.u-bordeaux1.fr/~delmoral/gips.html}{\tt Feynman-{K}ac formulae}.
\newblock \href{http://www.math.u-bordeaux1.fr/~delmoral/gips.html}{\tt Genealogical and interacting particle systems with applications}.
\newblock  Probability and its Applications (New York). (573p.) Springer-Verlag, New
  York (2004).

\bibitem{d-2013}
P. {Del Moral}.
\newblock \href{http://www.math.u-bordeaux1.fr/~pdelmora/Intro+Refs-Mean-Field-Simulation.pdf}{\tt Mean field simulation for Monte Carlo integration.} 
\newblock \href{http://www.crcpress.com/product/isbn/9781466504059}{\tt Chapman \& Hall.  Monographs on Statistics \& Applied Probability}  (2013).

\bibitem{dds1}
P. Del Moral, A. Doucet, and S. S. Singh. 
A backward interpretation of Feynman-Kac formulae. \emph{M2AN}, 44,
947--975 (2010).

\bibitem{dpr-2009}
P. {Del Moral}, F. Patras, and S. Rubenthaler. 
Coalescent tree based functional representations for some Feynman-Kac particle models.
{\em Ann. Appl. Probab.}, 19, 1--50 (2009).


\bibitem{dkp-14}
P. {Del Moral}, R. Kohn and F. Patras. On particle Gibbs Markov chain Monte Carlo models. \href{http://arxiv.org/abs/1404.5733}{arXiv:1404.5733 } (2014).

\bibitem{djlmp-2013}
P Del Moral, P. Jacob, A. Lee, L. Murray, and G.W. Peters 
Feynman-Kac particle integration with geometric interacting jumps, 
arXiv preprint arXiv:1211.7191 (2012). \emph{Stoch. Analysis Appl.}, 31, 830--871 (2013).

\bibitem {doucet2001}A. Doucet, J.F.G. de Freitas and N.J. Gordon Eds.
\emph{Sequential Monte Carlo Methods in Practice}. Springer-Verlag, New York (2001).

\bibitem{douc1}
R. Douc, and E. Moulines. Limit theorems for weighted samples with applications to sequential Monte Carlo methods. \emph{Ann. Statist.}, 36,  2344--2376 (2008).

\bibitem{douc}
R. Douc, A. Garivier, E. Moulines and J. Olsson. Sequential Monte Carlo smoothing for general state space Hidden Markov Models. \emph{Ann. Appl. Probab.}, 21, 2109--2145 (2011).

\bibitem {tony3}B. Jourdain, T. Leli\`evre and M. El Makrini. Diffusion Monte
Carlo method: numerical analysis in a simple case. \emph{Math. Model. Num.
Ana.}, 41, 189-213 (2007).


\bibitem {tony}T. Leli\`evre, M. Rousset and G. Stoltz. \emph{Free energy
computations: A mathematical perspective}. Imperial College Press (2010).

\bibitem{tony-rs}
T. Leli\`evre, M. Rousset and G. Stoltz. Computation of free energy differences through non-equilibrium stochastic dynamics: the reaction coordinate case. 
\emph{J. Comp. Phys.}, 222(2), 624-643, (2007).

\bibitem{ldm-2014}
F. Lindsten, R. Douc, and E. Moulines.
Uniform ergodicity of the Particle Gibbs sampler. \emph{Scand. J. Statist.} (to appear) (2015).


\bibitem{rennie}
B.C. Rennie, and A.J. Dobson. On Stirling numbers of the second kind. {\em J. Combinat. Theory}, 7, 116-121 (1969).


  \bibitem {rousset}
  M. Rousset. On the control of an interacting particle
approximation of Schr\"odinger ground states. \emph{SIAM J. Math. Anal.}, 
38, 824--844  (2006).

\end{thebibliography}
\end{document}